\documentclass[12pt]{amsart}
\usepackage{latexsym}
\usepackage{amssymb}
\usepackage{amscd}
\renewcommand\a{\alpha}
\renewcommand\b{\beta}
\newcommand\g{\gamma}
\renewcommand\d{\delta}
\newcommand\la{\lambda}

\newcommand\e{\eta}

\newcommand\s{\sigma}
\newcommand\x{\chi}

\newcommand\vf{\varphi}
\newcommand\p{\psi}
\renewcommand\t{\tau}
\renewcommand\r{\rho}

\newcommand\w{\omega}

\newcommand\vL{\varLambda}

\newcommand\ve{\varepsilon}

\newcommand{\PP}{\mathbb P}

\newcommand{\NN}{\mathbb N}

\newcommand\Fq{{\mathbf F}_q}

\newcommand\Ql{\bar{\mathbf Q}_l}

\newcommand\BF{\mathbf F}

\newcommand\BZ{\mathbf Z}

\newcommand\CB{\mathcal{B}}
\newcommand\ZC{\mathcal{C}}
\newcommand\CH{\mathcal{H}}
\newcommand\CI{\mathcal{I}}
\newcommand\CE{\mathcal{E}}

\newcommand\CM{\mathcal{M}}
\newcommand\CN{\mathcal{N}}
\newcommand\CO{\mathcal{O}}
\newcommand\CK{\mathcal{K}}
\newcommand\CP{\mathcal{P}}
\newcommand\CQ{\mathcal{Q}}

\newcommand\CF{\mathcal{F}}

\newcommand\CW{ \mathcal{W}}

\newcommand\Fg{\mathfrak g}

\newcommand\Fh{\mathfrak h}

\newcommand\Fl{\mathfrak l}

\newcommand\iv{^{-1}}

\newcommand\wt{\widetilde}
\newcommand\wg{^{\wedge}}

\newcommand\orth{^{\bot}}
\newcommand\ol{\overline}

\newcommand\lra{\leftrightarrow}

\newcommand\ssim{/\!\!\sim}

\newcommand\IC{\operatorname{IC}}
\newcommand\Ker{\operatorname{Ker}}
\newcommand\Hom{\operatorname{Hom}}
\newcommand\End{\operatorname{End}}

\newcommand\Res{\operatorname{Res}}

\newcommand\Tr{\operatorname{Tr}\,}

\newcommand\ad{\operatorname{ad}}

\newcommand\uni{_{\operatorname{uni}}}

\newcommand\lp{\operatorname{\!\langle\!}}
\newcommand\rp{\operatorname{\!\rangle\!}}

\newcommand\even{\operatorname{even}}
\newcommand\odd{\operatorname{odd}}
\newcommand\Spin{\operatorname{Spin}}

\newcommand\da{\dot a}

\newcommand{\isom}{\,\raise2pt\hbox{$\underrightarrow{\sim}$}\,}
\numberwithin{equation}{section}

\newtheorem{thm}{Theorem}[section]
\newtheorem{lem}[thm]{Lemma}
\newtheorem{cor}[thm]{Corollary}
\newtheorem{prop}[thm]{Proposition}

\def \para#1{\par\medskip\textbf{#1}
              \addtocounter{thm}{1}}

\def \remark#1{\par\medskip\noindent
                \textbf{Remark #1}
                \addtocounter{thm}{1}}

\begin{document}
\setlength{\baselineskip}{4.9mm}
\setlength{\abovedisplayskip}{4.5mm}
\setlength{\belowdisplayskip}{4.5mm}
\renewcommand{\theenumi}{\roman{enumi}}
\renewcommand{\labelenumi}{(\theenumi)}
\renewcommand{\thefootnote}{\fnsymbol{footnote}}
\renewcommand{\thefootnote}{\fnsymbol{footnote}}
\parindent=20pt
\medskip
\title{Generalized Green functions \\and unipotent classes
\\ for finite reductive groups, II }
\author{Toshiaki Shoji}
\par\medskip
\maketitle
\vspace{-0.5cm}
\par\vspace{1cm}
\begin{center}
Graduate School of Mathematics \\
Nagoya University  \\
Chikusa-ku, Nagoya, 464-8602,  Japan
\end{center}
\pagestyle{myheadings}
\markboth{SHOJI}{GENERALIZED GREEN FUNCTIONS}

\begin{abstract}
This paper is concerned with the problem of the determination of 
unknown scalars involved in the algorithm of computing the generalized Green
functions of reductive groups $G$ over a finite field.
In the previous paper, we have treated the case where $G = SL_n$.
In this paper, we determine the scalars in the case where $G$ 
is a classical group $Sp_{2n}$ or $SO_N$ for arbitrary characteristic.   
\end{abstract}

\bigskip
\medskip
\addtocounter{section}{-1}
\section{Introduction}
This paper is a sequel to [S2].
Our aim is to remove an ambiguity from the algorithm of
computing generalized Green functions of reductive groups
due to Lusztig.
Let $G$ be a connected reductive group defined over a finite 
field $\Fq$ with Frobenius map $F$.  Let $p$ be the
characteristic of $\Fq$.  
In [S2], we have treated the case where $G = SL_n$.  In this paper
we consider the case where $G = Sp_{2n}$ or $SO_N$ for arbitrary 
$p$.  The case where $G = \Spin_N$ will be treated in a separate paper.  
\par
In [S1] it was shown, in the case of $Sp_{2n}$ or $SO_N$ with 
$p \ne 2$, that there exists 
a representative in $C^F$ for each 
unipotent class $C$, called a distinguished element there 
(in this paper we call it a split element) which behaves well
with respect to the computation of Green functions.  Our result
in this paper shows that the split elements behave well 
for any type of generalized Green functions.  
We also show, in the case where $p = 2$,  that such a good representative
(called a split element) exists for $G = Sp_{2n}$ or $SO_{2n}$.
This was not known even for the case of Green functions.
\par
The main ingredient for the proof is a  variant of the restriction
theorem ([L1]) 
for the generalized Springer correspondence.  
The restriction theorem is a powerful tool for determining the 
generalized Springer correspondence, and it was used in [LS], [Sp2]
very effectively.  We extend this theorem so that it involves the 
information on the Frobenius action.  
In [S2], we have investigated the Frobenius action 
on the cohomology group $H^{a_0+r}_c(\CP_u, \dot\CE)$.
But this requires a precise information on the geometry of
$\CP_u$ related to the local system $\dot\CE$. 
In the case of classical groups, one can avoid to deal with 
$\CP_u$   
by considering the restriction theorem as above.

\section{A variant of the restriction theorem}
\para{1.1.}
We follow the notation in Section 1 in [S2].
In particular, $G$ is a connected reductive group over
a finite field $\Fq$, with Frobenius map $F$.  Let $k$ be 
an algebraic closure of $\Fq$ and $p$ the characteristic of $k$.
Let $P = LU_P$ be a parabolic subgroup of $G$, with a Levi 
subgroup $L$, and let $\CE$ be a cuspidal local system 
on a unipotent class $C$ in $L$.
As in (1.2.2) in [S2], one can define a perverse sheaf $K$
on $G$ associated to the triple $(L, C, \CE)$.
Then $K$ is a semisimple perverse sheaf with 
$\End K \simeq \Ql[\CW]$, where $\CW = N_G(L)/L$ is a Coxeter group.
Thus $K$ is decomposed as 
\begin{equation*}
\tag{1.1.1}
K = \bigoplus_{E \in \CW\wg}V_E\otimes K_E,
\end{equation*}
where
$K_E$ is a simple perverse sheaf on $G$ such that 
$V_E = \Hom(K_E, K)$ is an irreducible $\CW$-module corresponding
to $E \in \CW\wg$.
Put $d = \dim Z^0_L$, where $Z_L$ is the center of $L$. 
Let $G\uni$ be the unipotent variety of $G$, and $\CN_G$ the
set of all the pairs $(C',\CE')$, where $C'$ is a unipotent class 
in $G$ and $\CE'$ is a $G$-equivariant simple local system 
on $C'$.  Then it is known 
that $K[-d]|_{G\uni}$ is a semisimple perverse sheaf on $G\uni$, 
and it is decomposed as
\begin{equation*}
\tag{1.1.2}
K[-d]|_{G\uni} \simeq \bigoplus_{(C',\CE') \in \CN_G}
                  V_{(C',\CE')}\otimes \IC(\ol C', \CE')[\dim C'],
\end{equation*}
where $V_{(C',\CE')}$ is the multiplicity space for the simple
perverse sheaf $\IC(\ol C', \CE')[\dim C']$ on $G\uni$
 (cf. [S2, (1.2.4)]).  
Thus $K_E|_{G\uni}$ coincides with some $\IC(\ol C',\CE')$ up to
shift, and $V_{(C',\CE')}$ coincides with $V_E$.
It turns out that all the 
irreducible $\CW$-modules are realized as $V_{(C',\CE')}$ for 
some pair $(C',\CE')$.  Thus we have an injective map 
$\CW\wg \to \CN_G$ by $E = V_{(C',\CE')} \to (C', \CE')$, whose image
we denote by $\CN_G(C, \CE)$. 
Let $\CM_G$ be the set of triples $(L, C,\CE)$ up to $G$-conjugacy, 
where $L$ is a Levi subgroup of a parabolic subgroup of $G$ and $\CE$
is a cuspidal local system on a unipotent class $C$ of $L$.
The above injective maps form a bijection 
\begin{equation*}
\tag{1.1.3}
\coprod_{(L,C,\CE) \in \CM_G}(N_G(L)/L)\wg \to \CN_G
\end{equation*}
which is the so-called generalized Springer correspondence ([L1, 6.5]).
\para{1.2.}
Let $Q \supset P$ be a parabolic subgroup of $G$ with the Levi subgroup
$M$ such that $M \supset L$.  Then $\CW_1 = N_{M}(L)/L$ is 
in a natural way a subgroup of $\CW$.  Replacing $G$ by $M$, 
we have a subset $\CN_{M}(C,\CE)$ of $\CN_M$.  
For each $(C', \CE') \in \CN_G(C,\CE)$ (resp. 
$(C_1, \CE_1) \in \CN_{M}(C,\CE)$), we denote by 
$E$ (resp. $E_1$) the corresponding irreducible 
representation of $\CW$ (resp. ${\CW_1}$) under (1.1.3).
\par
Let $\pi_{Q} : Q \to M$ be the natural projection.  
Assume that $(C_1, \CE_1) \in \CN_{M}(C,\CE)$, and that
$(C',\CE') \in \CN_G$.  
We denote by
$f_{C_1, C'} : C_1U_Q \cap C' \to C_1$ the restriction  of 
$\pi_Q$. Then $\CF = R^{2d_{C_1,C'}}(f_{C_1,C'})_!\CE'$ is 
a semisimple $M$-equivariant local system on $C_1$,  
where $d_{C_1,C'} = (\dim C' - \dim C_1)/2$.
We define an integer $m_{\CE_1, \CE'}$ to be the 
multiplicity of $\CE_1$ in $\CF$. 
Lusztig proved the following restriction theorem on the 
generalized Springer correspondence.
\begin{thm}[{Lusztig [L1, Theorem 8.3]}] 
Under the above setting, $(C',\CE') \in \CN_G(C, \CE)$ if and only if  
$m_{\CE_1, \CE'} \ne 0$.  Moreover in that case we have 
\begin{equation*}
m_{\CE_1, \CE'} = \lp \Res E, E_1\rp_{\,\CW_1},
\end{equation*}
where $\lp\ ,\ \rp_{\,\CW_1}$ is the inner product of two representations of
$\CW_1$ (regarded as characters), and $\Res E$ is the restriction of 
$E$ on $\CW_1$.
\end{thm}
\para{1.4.}
Let $u \in C'$ and $v \in C_1$, and consider the component group
$A_G(u)$ and $A_M(v)$.  The set of $G$-equivariant simple local 
systems on $C'$ is in 1:1 correspondence with the set $A_G(u)\wg$ 
of irreducible characters of $A_G(u)$, 
and a similar fact holds also for $M$.    
As described in [LS], the integer $m_{\CE_1, \CE'}$ 
can be interpreted in terms of 
the representations of $A_G(u)$ and $A_M(v)$, which we explain below.
Let $\CF_v$ be the stalk of $\CF$ at $v \in C_1$.
Then we have
\begin{equation*}
\tag{1.4.1}
\CF_v \simeq H^{2d_{C_1, C'}}_c(C' \cap vU_Q, \CE').
\end{equation*}
Let $\pi : \wt C' = Z_G^0(u)\backslash G 
    \to  C', Z_G^0(u)g \mapsto g\iv ug 
$ 
be the finite covering of $C'$ with group $A_G(u)$.  
Let $X = (C' \cap vU_Q)\times_{C'} \wt C'$ be the fibre product
of $C' \cap vU_Q$ with $\wt C'$ over $C'$, and let 
$\wt \pi : X \to C' \cap vU_Q$ be the base change of $\pi$.
Then we have
\begin{equation*}
H_c^{2d_{C_1,C'}}(C' \cap vU_Q, \wt\pi_* \Ql) \simeq
       H_c^{2d_{C_1,C'}}(X, \Ql),
\end{equation*}
and $A_G(u)$ acts naturally on the right hand side.
Now $\wt\pi_*\Ql$ can be decomposed as 
$\wt\pi_*\Ql = \sum_{\r}V_{\r}\otimes \CE_{\r}$, where 
$\r$ runs over all the irreducible characters of $A_G(u)$. 
Here $\CE_{\r}$ is the $G$-equivariant simple local system on 
$C'$ corresponding to $\r$ and $V_{\r}$ is the corresponding irreducible
representation of $A_G(u)$. 
It follows that 
\begin{equation*}
H_c^{2d_{C_1,C'}}(C' \cap vU_Q, \CE_{\r}) \simeq 
     \bigl(H_c^{2d_{C_1,C'}}(X, \Ql)\otimes V_{{\r}^*}\bigr)^{A_G(u)},
\end{equation*}
where ${\r}^*$ is the dual representation of $\r$. 
On the other hand, the semisimple local system $\CF$ can be written as
$\CF = \sum_{\r_1}m_{\r_1}\CE_{\r_1}$, where $\CE_{\r_1}$ is the irreducible
local system on $C_1$ corresponding to $\r_1 \in A_M(v)\wg$ and $m_{\r_1}$ is
the multiplicity of $\CE_{\r_1}$ in $\CF$.    By taking the stalk at $v$, 
we have $\CF_v = \sum_{\r_1}m_{\r_1}(\CE_{\r_1})_v$.  
Here $(\CE_{\r_1})_v$ is an irreducible $A_M(v)$-module 
corresponding to $\r_1$.  
Note that if $\CE' = \CE_{\r}$, and $\CE_1 = \CE_{\r_1}$,  we have 
$m_{\CE_1,\CE'} = m_{\r_1}$.
Now $Z_M(v)$ acts on $C' \cap vU_Q$ by conjugation, and it induces 
an action of $A_M(v)$ on $H_c^{2d_{C_1,C'}}(C' \cap vU_Q, \CE')$.
We have
\begin{align*}
m_{\CE_1, \CE'} &= 
   \lp H_c^{2d_{C_1,C'}}(C' \cap vU_Q, \CE_{\r}), \r_1\rp_{A_M(v)} \\
   &= \lp H_c^{2d_{C_1,C'}}(X,\Ql)\otimes V_{{\r}^*})^{A_G(u)}, 
                        \r_1 \rp_{A_M(v)},
\end{align*} 
where $\lp\ ,\ \rp_{A_M(v)}$ denotes the inner product of characters of
$A_M(v)$.
\par
By Proposition 1.2 in [L1], it is known that $\dim X \le d_{C_1, C'}$.
Thus $H_c^{2d_{C_1,C'}}(X,\Ql)$ has a basis corresponding to the set of 
irreducible components of $X$ of dimension $d_{C_1,C'}$, and the action
of $A_G(u)$ on $H_C^{2d_{C_1,C'}}(X,\Ql)$ coincides with the
permutation action of $A_G(u)$ on those irreducible components of $X$. 
Since $\wt C' = Z^0_G(u)\backslash G$, we have 
\begin{align*}
X &= \{ (y, Z_G^0(u)g) \in (C' \cap vU_Q)\times \wt C' \mid 
          y = g\iv ug \} \\
  &= \{ Z_G^0(u)g \mid g\iv ug \in vU_Q\} \\
  &= Z_G^0(u)\backslash \{ g \in G \mid g\iv ug \in vU_Q\}.
\end{align*}
Put $Y_{u,v} = \{ g \in G \mid g\iv ug \in vU_Q\}$.
Then $Z_G(u) \times Z_M(v)$ acts on $Y_{u,v}$ by 
$(z, z') : g \mapsto zg{z'}\iv$ for $z \in Z_G(u), z' \in Z_M(v)$,
and the projection $Y_{u,v} \to X = Z_G^0(u)\backslash Y_{u,v}$ gives 
a bijection between the set of irreducible components of $X$ and
$Y_{u,v}$, which is compatible with the action of $A_G(u)$ and $ A_M(v)$.
Note that 
\begin{align*}
\dim Y_{u,v} &= \dim X + \dim Z_G^0(u) \\
             &= d_{C_1,C'} + \dim Z_G^0(u) \\
             &= (\dim Z_G(u) + \dim Z_M(v))/2 + \dim U_Q.
\end{align*}
Let $X_{u,v}$ be the set of irreducible components of $Y_{u,v}$ of
dimension $d_{C_1,C'} + \dim Z^0_G(u)$. 
It follows from the above discussion, we have 
\begin{cor}[Lusztig-Spaltenstein {[LS, 0.4, (4)]}]  
Let $\ve_{u,v}$ be the permutation representation of 
$A_G(u) \times A_M(v)$ on $X_{u,v}$.  Then we have
\begin{equation*}
\lp \Res E, E_1\rp_{\,\CW_1} = m_{\CE_1,\CE'} = \lp \ve_{u,v}, 
      \r\otimes \r_1^*\rp_{A_G(u)\times A_M(v)}.
\end{equation*}
\end{cor}
\para{1.6.}
We want to consider a variant of Corollary 1.5 which involves 
the Frobenius action.
Assume that $P$ is $F$-stable, and that the triple 
$(L, C, \CE) \in \CM_G$ is $F$-stable.  We choose 
$u_0 \in C^F$ and fix an 
isomorphism $\vf_0: F^*\CE \isom \CE$ so that the induced 
isomorphism $\CE_{u_0} \to \CE_{u_0}$ is of finite order. $\vf_0$ 
induces an isomorphism
$\vf: F^*K \isom K$.
For each pair $(C',\CE') \in \CN_G^F$, we choose $u \in {C'}^F$.
We fix an isomorphism 
$\p_{\CE'}: F^*\CE' \isom \CE'$ as follows; $F$ acts naturally on 
$A_G(u)$, and we consider the semidirect product 
$\wt A_G(u) = \lp\t\rp\ltimes A_G(u)$, where $\t$ is the restriction 
of $F$ on $A_G(u)$.  Since $(C',\CE')$ is $F$-stable, $\r$ is 
$F$-stable.  We choose an extension $\wt\r$ of $\r$ to $\wt A_G(u)$
and fix an isomorphism $\p_{\CE'}$ so that the induced isomorphism 
$\CE'_u \to \CE'_u$ corresponds to the action of $\t$ on $\wt\r$. 
Now $\p_{\CE'}$ 
induces an isomorphism 
$\wt\p_{\CE'}: F^*\IC(\ol C', \CE')[\dim C'] \isom 
         \IC(\ol C',\CE')[\dim C']$.
The isomorphism $\vf$ also induces an isomorphism 
$F^*K[-d]|_{G\uni} \isom K[-d]|_{G\uni}$, which we also denote by 
$\vf$.  Then under the decomposition of (1.1.2), $\vf$ induces 
an isomorphism 
\begin{equation*}
V_{(C'\CE')}\otimes F^*\IC(\ol C',\CE')[\dim C'] \isom 
 V_{(C'\CE')}\otimes \IC(\ol C',\CE')[\dim C']
\end{equation*}
for each pair
$(C',\CE') \in \CN_G^F$, and one can define a linear isomorphism 
$\s_{(C',\CE')}$ on $V_{(C',\CE')}$ such that 
this isomorphism can be written as
$\s_{(C',\CE')}\otimes \wt\p_{\CE'}$.
Now $F$ acts naturally on $\CW = N_G(L)/L$, and $\s_{(C',\CE')}$ 
becomes $\CW$-semilinear, namely we have
a relation $\s_{(C', \CE')}w = F(w)\s_{(C',\CE')}$ on $V_{(C',\CE')}$
for each $w \in \CW$. 
Replacing $G$ by $M$, we can define a $\CW_1$-semilinear map
$\s_{(C_1,\CE_1)}$ on $V_{(C_1,\CE_1)}$ for each pair 
$(C_1,\CE_1) \in \CN_M^F$.
The irreducible $\CW$-module $V_{(C',\CE')}$ can be written as 
a $\CW_1$-module 
\begin{equation*}
\tag{1.6.1}
V_{(C',\CE')} = \sum_{E_1 \in \CW_1\wg}M_{E_1}\otimes E_1,
\end{equation*}
where $M_{E_1}$ is the multiplicity space of the irreducible 
$\CW_1$-module $E_1$ and is realized as 
$M_{E_1} = \Hom_{\CW_1}(E_1, V_{(C',\CE')})$.
Suppose that $E_1 \simeq V_{(C_1,\CE_1)}$ under the generalized 
Springer correspondence for $M$.  If 
$E_1$ is $F$-stable, $(C_1, \CE_1) \in \CN_M^F$, and we have 
an isomorphism $\s_{(C_1,\CE_1)}$ on $E_1$. One can define a map 
$\s_{\CE_1,\CE'}: M_{E_1} \to M_{E_1}$ by 
$f \mapsto \s_{(C',\CE')}\circ f \circ \s_{(C_1,\CE_1)}\iv$.
The linear map $\s_{(C',\CE')}$ stabilizes the subspace 
$M_{E_1}\otimes E_1$ and we have 
\begin{equation*}
\tag{1.6.2}
\s_{(C',\CE')}|_{M_{E_1}\otimes E_1}
          = \s_{\CE_1,\CE'}\otimes \s_{(C_1,\CE_1)}.
\end{equation*}
\par
On the other hand, since $F(C') = C', F(C_1) = C_1$, 
the map $f_{C_1, C'}$ is $F$-equivariant.  Hence 
$\p_{\CE'}: F^*\CE' \isom \CE'$ induces an isomorphism 
$\p_{C_1,C'} : F^*\CF \isom \CF$, and so a linear isomorphism 
$\CF_v \to \CF_v$ which we denote by the same symbol 
$\p_{C_1,C'}$.
Now  the local system $\CF$ on $C_1$ corresponds to a representation 
$V$ of $A_M(v)$.  
$V$ can be decomposed as 
\begin{equation*}
V = \sum_{\r_1 \in A_M(v)\wg}M_{\r_1}\otimes \r_1,
\end{equation*}
where $M_{\r_1} = \Hom_{A_M(v)}(\r_1, V)$ is the multiplicity 
space of the irreducible 
representation $\r_1$.
$F$ acts on $A_M(v)$, and as in the case of $G$ 
we consider the semidirect product 
$\wt A_M(v) = \lp \t \rp \ltimes A_M(v)$, where $\t$ is the restriction 
of $F$ on $A_M(v)$.  
  For each $(C_1,\CE_1) \in \CN_M^F$, we 
fix an isomorphism $\p_{\CE_1}: F^*\CE_1 \isom \CE_1$ as in $G$ 
by using an extension $\wt\r_1$ of $\r_1$ to $\wt A_M(v)$.
Now $\p_{C_1,C'}$ stabilizes the subspace 
$M_{\r_1}\otimes \r_1$ for an $F$-stable $\r_1 \in A_M(v)\wg$,
and as in (1.6.2) 
one can define a linear map $\p_{\r_1,\r}$ on $M_{\r_1}$
such that 
\begin{equation*}
\tag{1.6.3}
\p_{C_1, C'}|_{M_{\r_1}\otimes\r_1} = 
       \p_{\r_1,\r}\otimes \p_{\CE_1}.  
\end{equation*}
\par
The following result gives an $F$-twisted version of the restriction 
theorem.  The proof is done by chasing the argument
in [L1].
\begin{prop} 
Under the notation as above, we have 
\begin{equation*}
\Tr(\s_{\CE_1,\CE'}, M_{E_1}) 
           = q^{-d_{C_1,C'} + \dim U_Q}\Tr(\p_{\r_1,\r}, M_{\r_1}).
\end{equation*}
\end{prop}
\para{1.8.}
Let $Y_{u,v}$ be as in 1.4. Assume that $Q$ is $F$-stable.  
Since $u,v$ are $F$-stable, $Y_{u,v}$ 
is $F$-stable, and so $F$ acts as a permutation on $X_{u,v}$. 
On the other hand, $F$ acts naturally on 
$A(u,v) = A_G(u) \times A_M(v)$, and we denote by $\wt A(u,v)$ the 
semidirect product group $\lp \t\rp \ltimes A(u,v)$. 
Then the permutation representation $\ve_{u,v}$ is extended to 
a representation of $\wt A(u,v)$, which we denote by $\wt\ve_{u,v}$. 
Now $\p_{\CE'}$ and $\p_{\CE_1}$ determines 
an extension of $\r\otimes \r_1^*$ to $\wt A(u,v)$, which we
 denote by 
$\wt{\r\otimes \r_1^*}$.  
By chasing the argument in 
1.4, we see that 
\begin{equation*}
\Tr(\p_{\r_1,\r}, M_{\r_1}) 
    = \lp\wt\ve_{u,v}, \wt{\r\otimes \r_1^*}\rp_{A(u,v)\t}, 
\end{equation*} 
where in general 
\begin{equation*}
\lp V_1, V_2\rp_{A(u,v)\t} = |A(u,v)|\iv\sum_{a \in A(u,v)}
         \Tr(a\t, V_1)\Tr((a\t)\iv, V_2)
\end{equation*}
for representations $V_1, V_2$ of $\wt A(u,v)$. 
Hence combined with Proposition 1.7, we have 
a variant of Corollary 1.5 involving the Frobenius action.
\begin{cor} 
Let the notations be as above. Then we have
\begin{equation*}
\Tr(\s_{\CE_1,\CE'}, M_{E_1}) = 
  q^{-d_{C_1,C'} + \dim U_Q}
    \lp\wt\ve_{u,v}, \wt{\r\otimes \r_1^*}\rp_{A(u,v)\t}.
\end{equation*}
\end{cor}
\para{1.10.}
We shall connect the above results to the discussion on 
generalized Green functions in [S2, Section 1].  
Take $j = (C',\CE') \in \CN_G^F$ and let  
$\p_0 = \p_{\CE'} : F^*\CE' \isom \CE'$ be defined in 1.6. 
$\p_0$ determines the $G^F$-invariant 
function $Y_j^0$ on the set $G^F\uni$ as in [S2, 1.3].
On the other hand, let $\wt\CW = \CW\rtimes \lp c\rp$ 
be the semidirect product, where $c$ is a Coxeter group 
automorphsim on $\CW$ induced from the action of $F$. 
In the decomposition in (1.1.1), one can define an isomorphism 
$\vf_E : F^*K_E \isom K_E$ so that the induced map 
$\s'_E : V_E \to V_E$ makes the irreducible $\CW$-module $V_E$ 
the preferred extension to $\wt\CW$ (cf. [L2, IV, (17.2)]). 
Put
\begin{align*}
a_0 &= -\dim Z_L^0 - \dim C', \\
r &= \dim G - \dim L + \dim\, (C \times Z_L^0).
\end{align*}
The we have 
\begin{equation*}
a_0 + r = (\dim G - \dim C') - (\dim L - \dim C).
\end{equation*}
We have $\CH^{a_0}(K_E)|_{C'} = \CE'$ 
and we define 
$\p: F^*\CE' \isom \CE'$ so that $q^{(a_0+r)/2}\p$ coincides with
the map defined by $\vf_E: F^*\CH^{a_0}(K_E) \isom K_E$. 
The function $Y_j$ is defined as the characteristic function 
of $\CE'$ through $\p$,  extended by 0 to the function on 
$G^F\uni$ (see [S2, 1.3 ]). 
Since $\CE'$ is a simple local system, there exists $\g \in \Ql^*$ 
such that $\p = \g \p_0$, and so $Y_j = \g Y_j^0$. 
Our main objective is the determination of this scalar $\g$.
Note that the determination of $\g$ is equivalent to the
determination of the map $\s_{(C',\CE')}$.  
In this paper, we determine $\g$ by investigating the map 
$\s_{(C',\CE')}$. 
The following 
fact is easily verified.
\begin{lem} 
Suppose that $q^{-(a_0+ r)/2}\s_{(C',\CE')}$ makes 
the $\CW$-module $V_{(C',\CE')}$ the preferred extension 
to $\wt\CW$.   Then we have $\g = 1$.
\end{lem}

\section{Unipotent classes of classical groups}
\para{2.1.}
Let $G$ be a connected classical group defined 
over $\Fq$.  
We consider the following type of groups $G$.
\par\medskip
(I) $G = Sp_{2n},$ \ $p \ne 2$,
\par\medskip
(II) $G = SO_{2n+1},$ \ $p \ne 2$,
\par\medskip
(III) $G = SO^{\pm}_{2n},$ \ $p \ne 2$.
\par\medskip
(IV) $G = Sp_{2n},$ \ $p = 2$, 
\par\medskip
(V) $G = SO^{\pm}_{2n},$ \ $p = 2$.
\par\medskip
These groups are realized as a group of transformations preserving
the various forms.  Let $V$ be a vector space over $k$ with 
$\dim V = N$. Assume that $p \ne 2$.
Then $Sp_{N}$ (resp. $O_{N}$) is the subgroup of 
$GL(V)$ leaving $f$ invariant, 
where $f$ is an alternating form (resp. a symmetric bilinear form ) 
on $V$ and $N = 2n$ in the case of $Sp$.  
$SO_N$ is the connected component of $O_N$, and $SO_{2n}^{\pm}$ 
corresponds to two $\Fq$-forms of $f$, one is split, the other 
is non-split.
\par 
Assume that $p = 2$.  Then $Sp_{2n}$ is the subgroup of 
$GL(V)$ with $N = 2n$ leaving an alternating form 
(= a symmetric bilinear form) $f$ invariant.
The quadratic form $Q$ on $V$ is defined by the property that 
the map 
$V \times V \to k, (x, y) \mapsto Q(x + y) - Q(x) - Q(y)$
gives rise to a non-singular bilinear form, which we may take
the alternating form $f$.
Let $O_{2n}$ be the subgroup of $GL(V)$ leaving $Q$ invariant. 
Then we have $O_{2n} \subset Sp_{2n}$, and let $SO_{2n}$ be 
the connected component of $O_{2n}$.   
It is known by [D] that there exists two $\Fq$-forms of $Q$
as follows.
We regard $Q$ as the quadratic form on $V_0 = \Fq^{2n}$.
Then there exists a basis of $V_0$ such that, 
for $x = (x_1, \dots, x_{2n}) \in V_0$ with respect to this 
basis, 
 $Q(x)$ can be expressed as 
\begin{align*}
\tag{2.1.1}
Q(x) &= x_1x_{n+1} + \cdots + x_{n}x_{2n},\\
\tag{2.1.2}
Q(x) &= x_1x_{n+1} + \cdots x_{n-1}x_{2n-1} + 
       \a x_n^2 + x_nx_{2n} + \a x_{2n}^2,
\end{align*}
where $\a \in \Fq$ is an element such that
$\a X^2 + X + \a$ is an irreducible polynomial
in $\Fq [X]$.
We denote by $O_{2n}^+$ (resp. $O_{2n}^-$) the group
$O_{2n}$ associated to the form in (2.1.1) (resp. (2.1.2)),
and let $SO_{2n}^{\pm}$ be the connected component of 
$O_{2n}^{\pm}$.
\para{2.2.}
We shall describe the unipotent classes in $G$.  As is well-known, 
in the case where $p \ne 2$, the unipotent classes of $G$ are 
described by unipotent classes in $GL(V)$ which are parametrized 
by partitions of $N$ through Jordan normal form. 
Let $\wt C_{\la}$ be the unipotent class in $GL(V)$ corresponding 
to a partition $\la$ of $N$.  We write $\la$ as 
$\la = (\la_1 \le \la_2 \le \cdots \le \la_r)$ with 
$\sum_{i = 1}^r\la_i = N$ or $\la = (1^{c_1}, 2^{c_2}, \dots)$, 
where $r = l(\la)$ is called the length of $\la$.
Assume that $G = Sp_{2n}$.  
Then $C_{\la} = \wt C_{\la} \cap G$ is non-empty if and only if 
$c_i$ is even for odd $i$, and in that
case $C_{\la}$ is a single conjugacy class in $G$.
While for $\wt G = O_N$, $C_{\la} = \wt C_{\la} \cap \wt G$ is 
non-empty if and only if $c_i$ is even for even $i$, and in that 
case $C_{\la}$ form a single class in $\wt G$.
Now $C_{\la}$ is already contained in $G = SO_N$, and so gives a
unipotent class in $G$  in almost
all cases.  The exceptions are the cases where $\la$ satisfies the 
condition ; 
$c_i = 0$ if $i$ is even, and $c_i$ is even for all odd $i$.  
In that case, 
$C_{\la}$ is divided into two classes $C'_{\la}$ and $C''_{\la}$ in 
$G$.
\para{2.3.}
In the case where $p = 2$, the parametrization of unipotent classes
is more complicated.  We shall describe it following Spaltenstein 
[Sp1, 2.6].
First assume that $G = Sp_{2n}$ with $p = 2$, and let $f$ be the 
associated alternating form.  Then the unipotent classes 
in $G$ are parametrized by a pair $(\la, \ve)$, where $\la$ is a
partition of $2n$ such that $c_i$ is even for odd $i$, and  
$\ve$ is an assignment $\ve : i \mapsto \ve_i \in \{ 0,1\}$ for 
even $i$ such that $c_i \ne 0$.  
Here $\ve_i = 1$ if $c_i$ is odd,  and $\ve_i = 0$ or 1
if $c_i$ is even. The correspondence with unipotent classes are 
given as follows.  Let $u$ be a unipotent element in $G$.  Then 
as an element in $GL(V)$, $u$ is parametrized by a partition 
$\la$ of $2n$, which satisfies a similar condition as in the case
of $p \ne 2$.
Now take even $i$ such that $c_i$ is even non-zero. 
We define a function $h_i$ on $\Ker (u-1)^i$ by 
$h_i(x) = f((u-1)^{i-1}x, x)$.  Then we put 
\begin{equation*}
\tag{2.3.1}
   \ve_i = \begin{cases}
               0  &\quad\text{ if } h_i \equiv 0, \\
               1  &\quad\text{ otherwise.}
            \end{cases}
\end{equation*}
The pair $(\la, \ve)$ is the one corresponding to the unipotent
class in $G$ containing $u$.  We denote by $C_{\la,\ve}$ the unipotent
class in $G$ corresponding to $(\la,\ve)$.
For a convenience sake, we extend $\ve$ to the function on $\NN$ by 
$\ve:  i \mapsto \ve_i$, where $\ve_i = \w$ for $i$ not
appeared above ($\w$ is a symbol not contained in $\{ 0,1\}$).
\par
Next assume that $G = SO_{2n}$ with $p = 2$.  We have 
$\wt G = O_{2n} \subset Sp_{2n}$.  Let 
$\wt C_{\la,\ve}$ be the unipotent class in $Sp_{2n}$ corresponding 
to $(\la, \ve)$.  Then $C_{\la,\ve} = \wt C_{\la,\ve} \cap \wt G$ is
a unipotent class in $\wt G$.  Thus unipotent classes in $\wt G$ 
are 1:1 correspondence with unipotent classes in $Sp_{2n}$.  
Now $C_{\la,\ve}$ is contained in $G$ if and only if 
$l(\la)$ is even.  Assume that $l(\la)$ is even.  Then 
$C_{\la,\ve}$ forms a single unipotent class in $G$ except for the case
where $c_i = 0$ for all odd $i$, and $\ve_i = 0$ for all even $i$
such that $c_i \ne 0$ (here $c_i$ is even for even $i$).
In the latter case, $C_{\la,\ve}$ splits into two classes 
$C'_{\la,\ve}$ and $C''_{\la,\ve}$ in $G$.
\para{2.4.}
Let $G$ be as in 2.1.  
For a convenience sake, we introduce a function 
$\ve$ on $\NN$ also 
in the case of $p \ne 2$.  Assume that $p \ne 2$.  
In the case of $G = Sp_{2n}$, we put
$\ve_i = 1$ if $i$ is even and $c_i \ne 0$, and put $\ve_i = \w$
otherwise.  In the case of $O_N$, we put $\ve_i = 1$ if $i$ is odd and
$c_i \ne 0$, and put $\ve_i = \w$ otherwise.   
Under this convention, we denote the class $C_{\la}$ in $Sp_{2n}$
or $SO_N$ by $C_{\la, \ve}$.
For $u \in G$, let $A_G(u)$ be the component group of $Z_G(u)$ 
as before.
In the case of $\wt G = O_{N}$, we also consider 
$A_{\wt G}(u)  = Z_{\wt G}(u)/Z^0_{\wt G}(u)$ for $u \in \wt G$.
Following [Sp1, 2.9], we shall describe the structure of 
$A_G(u)$ and $A_{\wt G}(u)$. 
\par
Assume that $G = Sp_{2n}$ with $p \ne 2$.  
Take $u \in C_{\la,\ve}$.  
We consider the generator $a_i$ corresponding to each $\la_i$.  
Then $A_G(u)$ is 
an abelian group generated by $a_i$ such that $\ve(\la_i) = 1$
under the condition that $a_i^2 = 1$ and that 
$a_i = a_j$ if $\la_i = \la_j$.
\par
Next assume that $\wt G = O_N$ with $p \ne 2$.  Take $u \in C_{\la,\ve}$. 
Then $A_{\wt G}(u)$ is an abelian group generated by $a_i$, exactly 
by the same condition as the case of $Sp_{2n}$.  
Now $A_G(u)$ is the subgroup of $A_{\wt G}(u)$ of index 2 generated 
by $a_ia_j$ for each $i \ne j$.
\par
Next assume that $G = Sp_{2n}$ with $p = 2$.
Take $u \in C_{\la,\ve}$. Again we consider the generators 
$a_i$ corresponding to $\la_i$.  Then $A_G(u)$ is an abelian group 
generated by $a_i$ such that $\ve(\la_i) \ne 0$ under the condition 
that $a_i^2 = 1$ and that $a_i = a_j$ if $\la_i = \la_j$ or
if $\la_i = \la_j+1$ or if $\la_i$ is even and $\la_i = \la_j +2$. 
\par
Finally assume that $\wt G = O_{2n}$ with $p = 2$, and $G = SO_{2n}$.
Take $u \in C_{\la,\ve}$.  Then $A_{\wt G}(u)$ is an abelian group 
generated by $a_i$, exactly by the same condition as the case of 
$Sp_{2n}$ with $p = 2$.  $A_G(u)$ is the subgroup of $A_{\wt G}(u)$
of index 2 generated by $a_ia_j$ for $i \ne j$.
\para{2.5.}
In what follows, we shall construct a normal form of unipotent 
elements in $G^F$.  As a preliminary for this, we consider the case
where $G = SO_{2n}$ with $p = 2$.  So assume given a vector space 
$V$ over $\Fq$ of dimension $N = 2n$ 
with a basis $e_1, \dots, e_N$, endowed with an alternating 
form $f$. We define an element $v \in GL(V)$ by $(v-1)e_j = e_{j-1}$
(with a convention $e_0 = 0$), and assume that $v$ leaves 
$f$ invariant.  We consider the following condition on $f$.
\begin{align*}
\tag{2.5.1}
&f(e_1, e_N) = 1, \\
\tag{2.5.2}
&f(e_i,e_N) = 0 \quad\text{ for } i = n+1, \dots, N, \\
\tag{2.5.3}
&f(e_i, e_k) + f(e_{i+1}, e_k) + f(e_i, e_{k+1}) = 0
\quad\text{ for } 0 \le i,k \le N-1.
\end{align*}
Note that (2.5.3) is equivalent to the condition that
$v$ leaves $f$ invariant.  
Also note that the conditions (2.5.1) $\sim$ (2.5.3) determines 
the alternating form $f$ invariant by $v$ uniquely. In fact, 
it follows from (2.5.2) and (2.5.3) that
\begin{equation*}
\tag{2.5.4}
f(e_i, e_j) = 0 \quad\text{ for $n+1 \le i,j \le N$}.
\end{equation*}
Also it follows from (2.5.1) and (2.5.3), we have
\begin{equation*}
\tag{2.5.5}
f(e_i,e_j) = \begin{cases}
                 0 &\quad\text{ if } i+j \le N, \\
                 1 &\quad\text{ if } i+j = N+1.
             \end{cases}
\end{equation*}
Hence it is enough to show that $f(e_i, e_j)$ is determined for
$1 \le i \le n$ and $n+1 \le j \le N$.
By (2.5.5) we have $f(e_n, e_{n+1}) = 1$.  Since 
$f(e_{n+1}, e_j) = 0$ for $j \ge n+1$, we have $f(e_n, e_j) = 1$
for $j \ge n+1$ by (2.5.3).
Then $f(e_i, e_j) = f(e_i, e_{j-1}) + f(e_{i+1}, e_{j-1})$ 
is determined for $1 \le i \le n$ by induction
on $j$ ($n+1 \le j \le N$).
\par
We consider a quadratic form $Q$ such that 
$Q(x+y) - Q(x) - Q(y) = f(x,y)$, which is left invariant by 
$v$. We have the following lemma. 
\begin{lem}
Let the notations be as above. Assume that 
$Q(e_N) = 0$.  Then $Q$ is determined uniquely, which is non-degenerate
of split type. 
\end{lem}
\begin{proof}
Since $Q$ is invariant by $v$, it is known by [Sp1, 6.10] that
$Q(e_i) = f_(e_i, e_{i+1})$ for $i = 1, \dots, N-1$.  Hence 
$Q$ is determined uniquely by $f$ and by the condition $Q(e_N) = 0$.
It is easy to see that this $Q$ actually gives rise to a quadratic form
invariant by $v$.   
In order to show that 
$Q$ is non-degenerate of split type, it is enough to see that there
exists a basis $e'_1, \dots, e_N'$ of $V$ satisfying the property 
\begin{align*}
\tag{2.6.1}
Q(e_i') &= 0 \quad\text{ for } i = 1, \dots, N, \\
f(e_i', e_j') &= \begin{cases}
                      1 &\quad\text{ if } i+j = N+1, \\
                      0 &\quad\text{ otherwise.}
                 \end{cases}
\end{align*}
We show (2.6.1). 
We consider the square matrix 
$A = (f(e_i, e_{N-j+1}))_{1\le i,j \le n}$ of degree $n$.
By (2.5.5) $A$ is a lower unitriangular matrix. 
For $k = 0, 1, \dots,$ we denote by $A_k$ the principal minor
matrix of $A$ of degree $2^k$.  We have $A_0 = (1)$.
We show that 
\par\medskip\noindent
(2.6.2) \ The matrix $A_k$ has the following property; 
for $k$ such that $2^{k+1} \le n$, we have 
\begin{equation*}
A_{k+1} = \begin{pmatrix}
             A_k & 0 \\
             A_k & A_k
          \end{pmatrix}.
\end{equation*}
If $2^k < n < 2^{k+1}$, $A$ is of the form
\begin{equation*}
A = \begin{pmatrix}
             A_k & 0 \\
             A'_k & A''_k
          \end{pmatrix},
\end{equation*}
where $A_k'$ is the minor matrix of $A_k$ of type 
$(n-2^k, 2^k)$ consisting of the first 
$(n-2^k)$-rows and all the columns, and $A_k''$ is the principal 
minor matrix of $A_k$ of degree $n-2^k$.  
\par\medskip
In fact, assume that $2^{k+1} \le n$.  By induction we may 
assume that 
(2.6.2) holds for $k-1$.  
Put $A = (a_{ij})$ with  $a_{ij} = f(e_i, e_{N-j+1})$.  Then 
by (2.5.3), we have 
\begin{equation*}
\tag{2.6.3}
a_{i,j} = a_{i-1,j} + a_{i-1, j-1}.
\end{equation*}
By induction, we see that the last row of $A_k$ 
is of the form $(1, \dots, 1)$.  Hence (2.6.3) implies 
that the $(2^k+1)$-th row of $A_{k+1}$ is of the form
$(1, 0, \dots, 0, 1, 0, \dots, 0)$ (1 appears in the first 
and the $(2^k+1)$-th coordinates), which coincides with the first row
of the matrix $(A_k, A_k)$.  Since the $(2^k+2)$-th row of $A_{k+1}$
is determined by $(2^k+1)$-th row by (2.6.3), and so on, 
we see that the minor matrix of $A_{k+1}$ of type $(2^k, 2^{k+1})$, 
consisting of the last $2^k$-rows and all the columns, coincides with
$(A_k, A_k)$.  (Note that since the last column of $A_k$ is of the form
${}^t(0, \dots, 0,1)$, the interaction between two $A_k$ does not occur 
in this computation).  Thus (2.6.2) holds for the case where 
$2^{k+1} \le n$.  The case where 
$2^k < n < 2^{k+1}$ is dealt similarly.  
\par
For $j$ such that $2^{a-1} < j \le 2^a$, we define a marked matrix
$A^{(j)}$ as follows. 
In the matrix $A_a$, the $(j,j)$ entry is contained in a minor matrix
$A_1 = \begin{pmatrix}
           1 & 0 \\
           1 & 1
        \end{pmatrix},$
where the $(j,j)$-entry corresponds to the $(1,1)$-entry 
(resp. $(2,2)$-entry) of $A_1$ if $j$ is odd (resp. even). 
We define a marked matrix $A_a^{(j)}$ by replacing 
the minor matrix $A_1$ in $A_a$ by $A_1^{\bullet}$, where
\begin{equation*}
A_1^{\bullet} = \begin{pmatrix}
                    1^{\bullet} &  0 \\
                    1^{\bullet} &  1^{\bullet}
                \end{pmatrix}
\quad\text{ or }\quad
A_1^{\bullet} = \begin{pmatrix}
                    1 &  0 \\
                    1 &  1^{\bullet}
                \end{pmatrix}
\end{equation*}
according as $(j,j)$ corresponds to $(1,1)$ or $(2,2)$ in $A_1$.   
In each of the matrices the marks $\bullet$ are attached to some 
entries in $A_a$.  For example, for $2< j \le 2^2$, $A_2^{(j)}$ is given as 
\begin{equation*}
\begin{pmatrix}
                1 & 0 & 0 & 0 \\
                1 & 1 & 0 & 0 \\
                1 & 0 & 1^{\bullet} & 0 \\
                1 & 1 & 1^{\bullet} & 1^{\bullet} 
\end{pmatrix}
\quad\text{ or }\quad
\begin{pmatrix}
                1 & 0 & 0 & 0 \\
                1 & 1 & 0 & 0 \\
                1 & 0 & 1 & 0 \\
                1 & 1 & 1 & 1^{\bullet} 
\end{pmatrix}
\end{equation*}
according to the case where $j = 3$ or $j = 4$.
For $k > a$, we define $A_k^{(j)}$ inductively as in (2.6.2) 
by replacing $A_k$ by $A_k^{(j)}$, starting from $A^{(j)}_a$.  
Then we define 
the matrix $A^{(j)}$ for $2^k \le n  < 2^{k+1}$ by replacing 
$A_k', A_k''$ by $(A_k^{(j)})', (A_k^{(j)})''$ which is defined similarly. 
\par
By a direct observation, we have
\par\medskip\noindent
(2.6.4) \ The matrix $A^{(j)}$ has the following properties.
\begin{enumerate}
\item 
In each row, the number of marked 1 is even except the $j$-th row, 
where the number is 1. 
\item
In each column containing the marked 1's, the entries except the marked 1
are all zero. 
\end{enumerate}
\par\medskip
We now define, for $j = 1, \dots, n $, the vector $e'_{N-j+1}$ by
\begin{equation*}
e_{N-j+1}' = \sum_ke_{N-k+1},
\end{equation*}
where the sum is taken over $1 \le k \le n$ such that
$k$-th column in $A^{(j)}$ contains a marked 1.
It follows from (2.6.4) that we have
\begin{equation*}
\tag{2.6.5}
f(e_i, e'_{N-j+1}) = \begin{cases}
                          1 &\quad\text{ if } i = j, \\
                          0 &\quad\text{ otherwise}
                     \end{cases}
\end{equation*}
for $1 \le i,j \le n$.
\par
We now consider the values of $Q$.  Since $Q$ satisfies 
the relation $Q(e_i) = f(e_i, e_{i+1})$,  it follows, by (2.5.4) and 
(2.5.5) together with our assumption that $Q(e_N) = 0$, that 
\begin{equation*}
\tag{2.6.6}
Q(e_i) = \begin{cases}
           0 &\quad\text{ if } i \ne n, \\
           1 &\quad\text{ if } i = n.         
         \end{cases}
\end{equation*}
Note that we have $e'_{n+1} = e_{n+1}$ by the previous computation.
Put $e'_n = e_n + e_{n+1}$. Then we have 
\begin{align*}
\tag{2.6.7}
Q(e_n') &= Q(e_n) + Q(e_{n+1}) + f(e_n, e_{n+1}) = 0, \\
f(e_n', e'_{n+1}) &= f(e_n + e_{n+1}, e_n) = 1.
\end{align*}
Now put $e_i' = e_i$ for $i = 1, \dots, n-1$.
Then by (2.5.4) and (2.5.5), together with (2.6.5) $\sim$ (2.6.7),
we see that the basis $\{ e_1', \dots, e_N'\}$ satisfies the 
relation (2.6.1). The lemma is proved. 
\end{proof}
\para{2.7.}
Let $G$ be as in 2.1.  We assume that $G^F$ is of split type. 
For each $F$-stable unipotent class $C$ in $G$, we shall construct 
a normal form $u$, called a split element, in $C^F$. The $G^F$-conjugacy 
class of $u$ is called the split class in $C^F$.
In the case where $p = 2$, we construct $u$ following 
[Sp1, II, 6.19]. 
First consider the case where $G = Sp_{2n}$ with $p = 2$.
Take a unipotent class $C_{\la,\ve}$ of $G$.
For each $j \ge 1$, put
\begin{equation*}
\tag{2.7.1}
n_j = \begin{cases}
          \la_j     &\quad\text{ if } \ve(\la_j) = 1, \\
          \la_j+1   &\quad\text{ if } \ve(\la_j) = \ve, \\
          \la_j+2   &\quad\text{ if }  \ve(\la_j) = 0.
      \end{cases}
\end{equation*}
We consider the vector space $V_j$ over $\Fq$ of dimension $n_j$
with basis $e^j_1, \dots, e^j_{n_j}$.
Assume that a non-degenerate alternating form ( $=$ a symmetric bilinear
form) $f_j$ on $V_j$ is given.  Let $H(V_j)$ be the subgroup of $GL(V_j)$ 
consisting of $g \in GL(V_j)$ which leaves the form $f_j$ invariant.
Put $H = Sp_{n_j}(k)$. Then $H(V_j)$ is regarded as a 
subgroup $H^F$ of $H$
under the natural $\Fq$-structure $F$ on $H$.
\par
One can construct $v_j \in H(V_j)$ such that 
$(v_j -1)e^j_i = e^j_{i-1}$ for $i = 1, \dots, n_j$
(under the convention that $e_0^j = 0$) with respect to the 
alternating form $f_j$ satisfying the property as given in 
2.5 (with $f = f_j, N = n_j$). 
\par
We also note that 
\medskip\par\noindent
(2.7.2) \  The image of 
$v_j \in Z_H(v_j)$ to $A_H(v_j)$ gives a generator $\bar a_j$ of
$A_H(v_j)$, where $A_H(v_j)$ is of order 1 or 2.
\par\medskip
For each $h = \la_j$, we shall construct a vector space $M_j$
over $\Fq$ with an alternating form $f_j$,  and $u_j \in H(M_j)$
as follows.
\par\medskip
(a) \ $\ve(h) = 1$.  
In this case, $h$ is even and
$c_h$ is odd or even.  We put $M_j = V_j$ and $u_j = v_j$.
Thus $u_j \in H(M_j)$ with respect to $f_j^0 = f_j$ on $M_j$.
Since $f_j^0(e_{n_j}^j, e_1^j) = 1$, we see that the function  
$x \mapsto f_j^0((u_j-1)^{h-1}x,x)$ is non-trivial on 
$\Ker (u_j-1)^h$. 
\par\medskip
(b) \ $\ve(h) = \w$. 
In this case $h$ is odd and $c_h$ is even.  Assume that 
$\la_j = \la_{j-1} = h$, and let $(V_j, f_j)$ and 
$(V_{j-1}, f_{j-1})$ be as before.  Note that 
$\dim V_j = \dim V_{j-1} = n_j = h+1$ by (2.7.1).
We define an alternating form $f'$ on $V_j\oplus V_{j-1}$
by the condition that $f'|_{V_j} = f_j$, $f'|_{V_{j-1}} = f_{j-1}$
and that $V_j \bot V_{j-1}$.
Let $L$ be a line in  $V_j\oplus V_{j-1}$ generated
by $ e_1^{j} + e_1^{j-1}$.
Then $L$ is an isotropic line by (2.5.5) and we put 
$M_j = L^{\bot}/L$. We have $\dim M_j = 2h$.
Since $v_j(e_1^j) = e_1^j, v_{j-1}(e_1^{j-1}) = e_1^{j-1}$, 
we see that $v_j + v_{j-1}$ fixes $L$, and so it induces 
a linear transformation on $M_j$, which we denote by $u_j$.
The form $f'$ induces an alternating form $f_j^0$.
We have $u_j \in H(M_j)$.
\par
By (2.5.1) and (2.5.5), $L\orth$ has a basis 
\begin{equation*}
e_n^j + e_n^{j-1}, e_{n-1}^j, e_{n-1}^{j-1}, \dots, 
      e_2^j, e_2^{j-1}, e_1^j + e_1^{j-1}.
\end{equation*}
Hence$L\orth/L$ has a basis 
\begin{equation*}
\bar e^j_n + \bar e_n^{j-1}, \bar e_{n-1}^j, \bar e_{n-1}^{j-1},
     \dots, \bar e_2^j, \bar e_2^{j-1}, \bar e_1^j = \bar e_1^{j-1},
\end{equation*}
where $\bar e_i^j, \bar e_i^{j-1}$ denote the image of 
$e_i^j, e_i^{j-1}$ on $(V_j\oplus V_{j-1})/L$.
\par\medskip
(c) \ $\ve(h) = 0$. 
In this case, $h$ is even and $c_h$ is even. 
Assume that $\la_j = \la_{j-1}$, and consider the vector 
spaces $V_j$ and $V_{j-1}$ as before. By (2.7.1), we have
$n_j = \dim V_j = \dim V_{j-1} = h + 2$.  We consider the 
alternating form $f'$ on $V_j \oplus V_{j-1}$ as before.
Let $N$ be the subspace of $V_{j} \oplus V_{j-1}$ spanned 
by $e_1^j + e_1^{j-1}$ and $e_2^j + e_2^{j-1}$.  Then 
$N$ is an isotropic subspace of $V_j\oplus V_{j-1}$ of 
dimension 2, and we put $M_j= N^{\bot }/N$. 
We have $\dim M_j = 2h$.  The alternating form $f'$ induces 
an alternating form $f^0_j$ on $M_j$. 
Now $v_j + v_{j-1}$ stabilized $N$, and so induces a linear
transformation on $M_j$ which we denote by $u_j$.
We see that $u_j \in H(M_j)$. 
\par
By (2.5.1) and (2.5.5), $N^{\bot}$ 
has a basis
\begin{equation*}
e_n^j + e_n^{j-1}, e_{n-1}^j + e_{n-1}^{j-1}, 
    e_{n-2}^j, e_{n-2}^{j-1}, \dots, e_1^j, e_1^{j-1}. 
\end{equation*}
Hence $N^{\bot}/N$ has a basis 
\begin{equation*}
\bar e_n^j + \bar e_n^{j-1}, \bar e_{n-1}^j + \bar e_{n-1}^{j-1}, 
    \bar e_{n-2}^j, \bar e_{n-2}^{j-1}, \dots, 
\bar e_3^{j}, \bar e_3^{j-1}, \bar e_2^j = \bar e_2^{j-1},
\bar e_1^j = \bar e_1^{j-1}, 
\end{equation*}
where $\bar e_i^j, \bar e_i^{j-1}$ denotes the image 
of $e_i^j$ on $(V_j \oplus V_{j-1})/N$. 
The action of $u_j$ on this basis is easily described, and 
by using the formulas in 2.5, one can check that 
$f^0_j((u_j-1)^{h-1}x, x) = 0$ for all $x \in \Ker(u_j-1)^h$.
For example, 
\begin{equation*}
f^0_j(\bar e_n^j + \bar e_n^{j-1}, \bar e_3^j + \bar e_3^{j-1})
  = 2f_j(e_n^j, e_3^j) = 0,
\end{equation*}
and the other cases are dealt similarly.
\par\medskip
We now define a vector space $\bar V$ as $\bar V = \bigoplus_jM_j$ 
so that $\dim \bar V = 2n$, and let $f = \bigoplus_j f^0_j$ be 
the alternating form obtained 
from $f_j^0$.  Put $\bar u = \prod_j u_j \in H(\bar V)$.
It follows from the previous construction,  we have
\par\medskip\noindent
(2.7.3) \ $H(\bar V)$ can be identified with $G^F$, and under this
isomorphism, the element $u \in G^F$ corresponding to 
$\bar u$ gives an element in $C_{\la,\ve}^F$. 
We call $u$ a split element in $C_{\la,\CE}^F$.
\par\medskip
The structure of the group $A_G(u)$ is also described 
as follows (cf. [Sp1, II, 6.19])
Take $\la_j$ such that $\ve(\la_j) \ne 0$.  
We denote by  $\bar a_j$ an automorphism of  
$\bigoplus_{k \ge 1} V_k$ defined by
$v \mapsto v_j(v)$ for $v \in V_j$, and $v \mapsto v$ for $v \in V_k$
such that $k \ne j$.  
Then $\bar a_j$ induces an automorphism on $\bar V$ commuting with 
$\bar u$, which we denote also by $\bar a_j$.  
It is checked that $\bar a_j \in H(\bar V)$, and so this gives an 
element of $Z_G(u)$.  The image of $\bar a_j$ on $A_G(u)$ 
coincides with the generator $a_j$ stated in 1.4. 
In particular, we see that $F$ acts trivially on $A_G(u)$.
This implies, since $A_G(u)$ is abelian, that
\par\medskip\noindent
(2.7.4) \ For any $u' \in C_{\la,\ve}^F$, $F$ acts trivially on 
$A_G(u')$. 
\para{2.8.}
Next we consider the case where $G = SO_{2n}$ with $p = 2$.
Take a unipotent class $C_{\la,\ve}$ of $G$, and let 
$n_j$ be as in (2.7.1).
We consider the vector space $V_j$ over $\Fq$ with basis 
$e_1^j, \dots, e_{n_j}^j$ and a unipotent element 
$v_j \in GL(V_j)$ as in 2.7.  By Lemma 2.6, one can construct 
a split quadratic form $Q_j$ on $V_j$ which is invariant by $v_j$.
Let $H(V_j)$ be the subgroup of $GL(V_j)$ consisting of $g$ which
leaves $Q_j$ invariant.  Let $\wt H = O_{n_j}(k)$.  Since 
$Q_j$ is of split type, $H(V_j)$ can be identified with the subgroup 
$\wt H^F$ of $\wt H$, where $F$ is a split Frobenius map.
Then by a similar argument as in 2.7, we obtain 
$u_j \in H(M_j)$ for each case (a), (b) or (c), where $H(M_j)$ is the 
group of invariants with respect to the induced quadratic from 
$Q_j^0$. Note that the explicit computation in the proof of Lemma 2.6
shows that $Q_j^0$ is of split type.  As in 2.7, we define a vector 
space $\bar V = \bigoplus_j M_j$ and $Q = \bigoplus_j Q_j^0$, and put 
$\bar u = \prod_j u_j \in H(\bar V)$.  We have
\par\medskip\noindent
(2.8.1) \ Let $\wt G = O^+_{2n}$.  Then 
$H(\bar V)$ can be identified with $\wt G^F$ with split Frobenius 
map $F$,  
and under this isomorphism, the element $u \in \wt G^F$ corresponding to
$\bar u$ gives an element in $\wt C_{\la,\ve}^F$. 
In the ordinary case $u \in C_{\la, \ve}^F$. In the exceptional case
we have $u \in (C_{\la,\ve}')^F$ and $\s(u) \in (C_{\la,\ve}'')^F$,
where $\s$ is the graph automorphism on $SO_{2n}$. 
We call $u$ and $\s(u)$ the split elements in $\wt C_{\la}^F$. 
\para{2.9.}
Next we consider the case where $G = Sp_{2n}$ or $SO_N$ with $p \ne 2$.
We assume that $F$ is a split Frobenius map.
Let $\Fg$ be the Lie algebra of $G$. Since the unipotent
classes in $G$ are in bijection with the nilpotent orbits in $\Fg$
with $\Fq$-structure, we consider the normal form of nilpotent orbits
instead of unipotent classes. 
Let $\CO_{\la,\ve}$ be the nilpotent orbit in $\Fg$ corresponding to 
the unipotent class $C_{\la,\ve}$ in $G$. 
For each $\la_j$, we construct a vector space $M_j$ over $\Fq$ 
and a nilpotent 
transformation $X_j$ on $M_j$ as follows.  
\par\medskip
(a) $\ve(\la_j) = 1$.  We consider a vector space
$M_j$ of dimension $h = \la_j$ with basis $e^j_1, \dots, e^j_h$.
We define a non-degenerate alternating form (resp. a symmetric bilinear
form) $f_j$ on $M_j$ in the case where $G = Sp_{2n}$ (resp. $SO_N$) 
by  
\begin{equation*}
\tag{2.9.1}
f_j(e^j_{h-i+1}, e^j_i) = (-1)^{\d_j - i} \quad\text{ for }
i = 1, \dots, h,
\end{equation*}
where 
\begin{equation*}
\d_j = \begin{cases}
         \la_j/2 + j &\quad\text{ if } G = Sp_{2n}, \\
         (\la_j-1)/2 + j &\quad\text{ if } G = SO_N.
       \end{cases}
\end{equation*}
We put the value of $f_j$ zero  for any 
other pair of the basis.
We define a nilpotent transformation $X_j$ on 
$M_j$ by $X_j(e^j_i) = e^j_{i-1}$ 
(under the convention that $e_0 = 0$).  Then $X_j \in \Fh(M_j)$,
where $\Fh(M_j)$ is the subalgebra of 
$\Fg\Fl(M_j)$ consisting of
$X$ such that $f_j(Xx, y) + f_j(x, Xy) = 0$ for 
$x,y \in M_j$. 
\par\medskip
(b) \ $\ve(\la_j) = \w$.  In this case $c_h$ is even for $\la_j = h$.
We assume that $\la_j = \la_{j-1}$.
We consider a vector space $M_j$ of dimension $2h = 2\la_j$ with basis 
$e^j_1, \dots, e_{h}^j, e_1^{j-1}, \dots, e_h^{j-1}$.
We define an alternating form (resp. a symmetric bilinear form) $f_j$ on 
$M_j$  in the case where $G = Sp_{2n}$ (resp. $G = SO_N$) by 
\begin{equation*}
\tag{2.9.2}
f_j(e^j_{h-i+1}, e^{j-1}_i) = \ve f_j(e^{j-1}_i, e^j_{h-i+1})
                              = (-1)^{i-1} \quad\text{ for }
                        i = 1, \dots, h,
\end{equation*}  
where $\ve = -1$ (resp. $\ve = 1$) if $G = Sp_{2n}$ (resp. $G = SO_N$).
We put the values of $f_j$ zero for any other pair of the basis.
We define a nilpotent transformation $X_j$ on $M_j$ by 
$X_je_i^j = e^j_{i-1}, X_je^{j-1}_i = e^{j-1}_{i-1}$ (we put 
$e_0^j = e_0^{j-1} = 0$ as before).  Then $X_j \in \Fh(M_j)$. 
\par\medskip
We define a vector space $\bar V$ by $\bar V = \bigoplus_j M_j$
so that $\dim \bar V = N$, and let $f = \sum_j f_j$ be an 
alternating form (resp. a symmetric bilinear form) on $\bar V$ obtained
from $f_j$.  Put $\bar X = \bigoplus_j X_j \in \Fh(\bar V)$.  
Then it is known by [SS], 
\par\medskip\noindent
(2.9.3) \  $\Fh(\bar V)$ can be identified with $\Fg^F$. Under this 
correspondence $\bar X$ gives an element $X \in \CO_{\la,\ve}^F$,
which we call a split element in $\CO_{\la,\ve}^F$.  $X$ also determines
the $G^F$-class in $C_{\la,\ve}^F$, which we call the split class 
in $C^F_{\la,\ve}$.
\para{2.10.}
We consider the case where $G = SO_{2n}$ with non-split Frobenius $F$
(for arbitrary $p$).
Let $F_0$ be the split Frobenius map.  Then one can write $F = F_0\s$
with the graph automorphism $\s$.  We may choose 
$s \in \wt G\  \backslash\  G$ such that $\s = \ad s$, and fix it for 
all.
Let $u \in G^{F_0}$ be a split element in  $C_{\la,\ve}$.
Since $A_{\wt G}(u) \ne A_G(u)$, there exists 
$a \in A_{\wt G}(u)\ \backslash \ A_G(u)$.  Let $\da \in Z_{\wt G}(u)$
be a representative of $a$.  Since $[\wt G : G] = 2$, 
there exists $g \in G$ such that $\da = gs$.  It follows that 
${}^{gF}u = u$.  Now there exists $\a \in G$ such that 
$\a\iv F(\a) = g$, and we have $u' = \a u\a\iv \in C_{\la, \ve}^F$.
It is easy to see that the $G^{F}$-conjugacy class of $u'$ is uniquely 
determined by $\da \in Z_{\wt G}(u)$.  For the exceptional case, 
we have $u' \in (C_{\la,\ve}')^F$ and $\s(u') \in (C_{\la,\ve}'')^F$.
In what follows, we fix a split element $u' \in C_{\la,\ve}^F$ as
follows.
\par\medskip\noindent
(2.10.1) \ Let $G = SO^-_{2n}$.  
Let $u$ be the split element in $C_{\la,\ve}^{F_0}$ or 
$(C_{\la,\ve}')^{F_0}$. We choose 
$a = a_i \in A_{\wt G}(u)$ such that $\ve(\la_i) = 1$ (resp. 
$\ve(\la_i) \ne 0$) in the case where $p \ne 2$ (resp. $p = 2$)
and that $\la_i$  is minimal under this condition.  Take a 
representative $\da_i \in Z_{\wt G}(u)$ as in 2.7, and define 
$u' \in C_{\la,\ve}^F$ by using $\da_i$.
(In the exceptional case, define $u' \in (C_{\la,\ve}')^F$ and 
$\s(u') \in (C_{\la,\ve}'')^F$.)  We call $u'$ the split element
in $C_{\la,\ve}^F$.
\par\medskip
Under the notation above, $aF_0 = gF$ acts trivially on 
$A_G(u)$ by (2.7.4).  It follows that $F$ acts trivially on 
$A_G(u')$.  Since $A_G(u')$ is abelian, we have
\par\medskip\noindent
(2.10.2)\ The statement (2.7.4) holds also for the case where 
$F$ is of non-split type. 
\remark{2.11.}
The definition of the split elements for $Sp_{2n}$ or $SO_N$ 
(for $p \ne 2$) in 2.9, 2.10 
is essentially the same as the one used in [S1, 3.3, 3.7]
(where it is called distinguished elements). 
Note that in the case where $q \equiv 1 \pmod 4$, $\d_j$ can
be removed in the formula (2.9.1) by a suitable base change. 
Also note that the definition of $f_j$ involves the case where
$\ve(h) = 1$ and $c_h$ is even, which is necessary for later
discussions, though these cases are ignored in [S1]. 
\par
In the case of non-split groups with $p \ne 2$, 
our definition of split elements
is not the same as in [S1, 3.7], where it is defined by using 
$a = a_i$ corresponding to $\la_i$ of maximal length 
instead of minimal length.  This is not essential, but 
the definition here is more convenient since it produces 
preferred extensions of $\CW$-modules as will be seen in Theorem 4.3. 
(The elements defined in [S1, 3.7] do not necessarily 
produce them).
\par
\newpage
\section{Generalized Springer correspondence}
\para{3.1.}
We review here the generalized Springer correspondence for classical
groups following [L1] and [LS].  In what follows, we denote by 
$W_n$ the Weyl group of type $C_n$, and by $W_n'$ the Weyl group of 
type $D_n$.  Throughout the whole cases, for a given 
$(L,C, \CE) \in \CM_G$, the cuspidal pair $(C, \CE)$ is uniquely determined.
So, we just describe $L$ which has a cuspidal pair. 
\par\medskip
(a) \ 
Let $G = Sp_{2n}$ with $p \ne 2$.  Then $(L,C,\CE) \in \CM_G$ if 
and only if $L$ is of type $C_m$ for some $m$ of the form 
$m = \frac{1}{2}d(d-1)$ with $d \ge 1$.  
We have $N_G(L)/L \simeq W_{n-\frac{1}{2}d(d-1)}$.  
Since the set $\{ d(d-1) \mid d \ge 1\}$ coincides with the set
$\{ d(d-1) \mid d \in \BZ, d : \text{odd}\}$, the generalized Springer
correspondence (1.1.3) is given by a bijection 
\begin{equation*}
\CN_G \lra \coprod_{\substack{d \in \BZ \\ 
                     d\text{ odd}}}
         (W_{n-\frac12 d(d-1)})\wg.
\end{equation*} 
\par
(b) \ Let $G = SO_N$ with $p \ne 2$.  
Then $(L,C,\CE) \in \CM_G$ if and only if $L$ is of type 
$B_m$ (resp. $D_m$) for some $m$ such that $m = \frac 12 (d^2-1)$ 
(resp. $m = \frac 12 d^2$) 
and that $d \equiv N \pmod 2$ in the case where $N$ is odd 
(resp. $N$ is even) with $m \ge 0$.  
We have 
$N_G(L)/L \simeq W_{(N - d^2)/2}$ if $m \ge 1$, and 
$N_G(L)/L \simeq W_n$ (resp. $W_n'$) if $m = 0$, namely if 
$L$ is a maximal torus $T$, in the case where $N$ is odd 
(resp. $N$ is even). Thus the generalized Springer correspondence 
is given by a bijection
\begin{align*}
\CN_G &\lra \coprod_{\substack{d \ge 1 \\ 
       d\text{ odd} }}(W_{(N - d^2)/2})\wg
       \quad\text{ ($N$ odd)}, \\ 
\CN_G &\lra W_n' \coprod\bigl( \coprod_{\substack{ d > 0 \\
           d\text{ even} }}(W_{(N -d^2)/2})\wg\bigr)
       \quad\text{ ($N$ even)}. 
\end{align*}
\par
(c) \ Let $G = Sp_{2n}$ with $p = 2$.
Then $(L, C,\CE) \in \CM_G$ if and only if $L$ is of type $C_m$
for some $m$ of the form $m = d(d-1)$ with $d \ge 1$.
We have $N_G(L)/L \simeq W_{n - d(d-1)}$.  
Hence as in the case (a), the generalized Springer correspondence 
is given by 
\begin{equation*}
\CN_G \lra \coprod_{\substack{ d \in \BZ \\ d\text{ odd}}}
               (W_{n-d(d-1)})\wg.
\end{equation*}
\par
(d) \ Let $G = SO_{2n}$ with $p = 2$.
Then $(L,C,\CE) \in \CM_G$ if and only if $L$ is of type $D_m$
for some $m$ of the form $m = d^2$ with $d \ge 0$, even.
We have $N_G(L)/L \simeq W_{n-d^2}$ if $d \ge 1$ and 
$N_G(L)/L \simeq W_n'$ if $d = 0$, namely if $L$ is a maximal torus
of $G$.  Hence the generalized Springer correspondence is given by
\begin{equation*}
\CN_G \lra W_n'\coprod \bigl(\coprod_{\substack{
        d > 0 \\ d\text{ even}}}(W_{n-d^2})\wg\bigr).
\end{equation*}
\para{3.2.}
The generalized Springer correspondence for classical groups 
is described in terms of symbols.
We review the notion of symbols following [L1], [LS].
Let $r, s \in \BZ_{\ge 1}$, $d \in \BZ$.  For each integer 
$n \ge 1$ let $\wt X_{n,d}^{r,s}$ be the set of all ordered
pairs $(A, B)$ where $A = \{ a_1, \dots, a_{m+d}\}$, 
$B = \{ b_1, \dots, b_m\}$ (for some $m$) are subsets of 
$\BZ_{\ge 0}$ satisfying the following conditions.
\begin{align*}
a_i - a_{i-1} &\ge r+s \quad (1 < i \le m+d), \\
\tag{3.2.1}
b_i - b_{i-1} &\ge r+s \quad (1 < i \le m), \\
b_1 &\ge s, \\
\sum a_i + \sum b_i &= n + (r+s)(m + [d/2])(m + d - [d/2])
       -r(m+ [d/2]). 
\end{align*}
(In the case where $r+s = 0$, $A$ or $B$ contains elements
with multiplicities.  In that case, we regard it as a sequence of 
integers. )
\par
Note that if $r = s = 1$ and $d$ is odd the fourth condition 
is written as 
\begin{equation*}
\tag{3.2.2}
\sum a_i + \sum b_i = n + \frac 12 (2m+d)(2m+d -1),
\end{equation*}
and if $r = 2, s = 0$ it is written as 
\begin{equation*}
\tag{3.2.3}
\sum a_i + \sum b_i = n + \frac 12 ((2m+ d -1)^2 -1).
\end{equation*}
\par
Let $X_{n,d}^{r,s}$ be the set of equivalence classes on 
$\wt X_{n,d}^{r,s}$ for the equivalence relation generated by
\begin{equation*}
(A,B) \sim \bigl(\{ 0\} \cup (A + r+s), \{ s\}\cup (B + r+s)\bigr),
\end{equation*}
where $A + r + s$ denotes the set 
$\{ a_1+(r+s), \dots, a_{m+d}+ (r+s)\}$ and so on for $B$.
We put
\begin{equation*}
X_n^{r,s} = \coprod_{d\text{ odd}}X_{n,d}^{r,s}.
\end{equation*}
An element in $X_{n,d}^{r,s}$ is called an $(r,s)$-symbol of 
rank $n$ and defect $d$, which we also denote by $(A,B)$.
\par
We consider the special case where $s = 0$.  In that case,
$(A,B) \mapsto (B,A)$ defines a bijection from 
$X_{n, d}^{r,0}$ to $X_{n,-d}^{r,0}$ and so induces an involution 
of each of the following sets,
\begin{equation*}
X_{n,\text{ even}}^{r,0} = \coprod_{d \text{ even}}X_{n,d}^{r,0}, 
   \quad X_{n,\text{ odd}}^{r,0} = 
        \coprod_{d \text{ odd}}X_{n,d}^{r,0}.
\end{equation*}
Let $Y^r_{n,\text{ odd}}$ (resp. $Y^r_{n,\text{ even}}$) 
be the set obtained as the quotient of $X_{n,\text{ odd}}^{r,0}$
(resp. $X_{n,\text{ even}}^{r,0}$) by this involution, with the
convention that the symbol invariant by the involution, i.e., 
the symbol $(A,A)$ which we call the degenerate symbol, is counted twice. 
For $d \ge 0$, the image of $X_{n,d}^{r,0}$ in $Y^r_{n,\text{ odd}}$
or $Y^r_{n,\text{ even}}$ is denoted by $Y_{n,d}^r$.
One can regard the element in $Y_{n,d}^r$ as a symbol $(A, B)$ in 
$X_{n,d}^{r,0}$ considered as an unordered pair.
\para{3.3.}
The set $W_n\wg$ is parametrized by ordered pairs of partitions
$(\a, \b)$ such that $|\a| + |\b| = n$.  For a fixed $d \in \BZ$, 
one can express the partitions $\a, \b$ as 
$\a: \a_1 \le \a_2 \le \cdots \le \a_{m+d}$, 
$\b: \b_1 \le \b_2 \le \cdots \le \b_m$ for a suitable $m$, by allowing 
0 in the entries.
Then $(A, B) \in \wt X_{n,d}^{0,0}$ for 
$A = \{ \a_1, \dots, \a_{m+d}\}$, $B = \{ \b_1, \dots, \b_m\}$,
and this induces a well-defined bijection between $W_n\wg$ and 
$X_{n,d}^{0,0}$. 
The same map induces a bijection between $W_n\wg$ and $Y_{n,d}^0$ if 
$d \ge 1$. 
On the other hand, the set $(W_n')\wg$ is parametrized by unordered
pairs of partitions $(\a, \b)$ such that $|\a| + |\b| = n$, under the
convention that $(\a,\a)$ is counted twice.  Thus in a similar way 
as above, we have a natural bijection between $(W_n')\wg$ and $Y_{n,0}^0$.
\par
For a given $r,s, d$, we define a symbol $\vL_d^{r,s} = (A, B)$ as
follows.
\begin{equation*}
     \begin{cases}
A = \{ 0, (r+s), \dots, (d-1)(r+s)\}, \quad B = \emptyset 
    &\quad\text{ if } d > 0, \\
A = \emptyset, \quad B = \{ s, s+(r+s), \dots, s+ (-d-1)(r+s)\}
    &\quad\text{ if } d < 0, \\
A = \emptyset, \quad B = \emptyset
    &\quad\text{ if } d = 0.
     \end{cases}
\end{equation*}
It is easy to see that $\vL_d^{r,s} \in X_{n_0,d}^{r,s}$ with 
$n_0 = (r+s)[d/2](d - [d/2]) - s[d/2]$, and that 
the set $X_{n_0,d}^{r,s}$ consists of a unique element $\vL_d^{r,s}$. 
In the case where $s = 0$, let $\vL_d^r$ be the image of 
$\vL_{d}^{r,0}$ under the map $X_{n_0,d}^{r,0} \to Y_{n_0,d}^r$.
For each $d \in \BZ$, one can define a map 
\begin{equation*}
\tag{3.3.1}
X_{n-n_0,d}^{0,0} \to X_{n,d}^{r,s}, \quad \vL \mapsto \vL + \vL_d^{r,s},
\end{equation*}
which gives a bijection $X_{n-n_0,d}^{r,s} \simeq X_{n,d}^{r,s}$.
(Note, in general, that the sum of two symbols $\vL, \vL'$ with the 
same defect is defined by choosing representatives 
$\vL = (A, B), \vL' = (A', B')$ such that  
$|A| + |B| = |A'| + |B'|$, namely of the same shape, 
and then by adding entry-wise.)
Similarly for each $d \ge 0$, one can define a bijection
\begin{equation*}
\tag{3.3.2}
Y_{n-n_0,d}^0 \to Y_{n,d}^r, \quad \vL \mapsto \vL + \vL_d^r.
\end{equation*}
\par
Combining (3.3.1) with the bijection 
$W_{n-n_0}\wg \simeq X_{n-n_0,d}^{0,0}$ above, we have a bijection
\begin{equation*}
\tag{3.3.3}
W_{n-n_0}\wg \to X_{n,d}^{r,s}.
\end{equation*}
Similarly, combining (3.3.2) with the bijections 
$W_{n-n_0}\wg \simeq Y_{n-n_0,d}^0$ for $d > 0$ and  
$(W_{n-n_0}')\wg \simeq Y_{n-n_0,0}^0$,  we have bijections
\begin{align*}
\tag{3.3.4}
W_{n-n_0}\wg &\to Y_{n,d}^r \quad (d > 0), \\ 
(W_{n-n_0}')\wg &\to Y_{n, d}^r \quad (d = 0). 
\end{align*}
\para{3.4.}
A symbol $(A,B) \in X_{n,d}^{r,s}$ is said to be distinguished 
if $d = 0$ or 1, and  
\begin{equation*}
\begin{aligned}
&a_1 \le b_1 \le a_2 \le b_2 \le \cdots \le a_m \le b_m 
         \quad &(d = 0),  \\   
&a_1 \le b_1 \le a_2 \le b_2 \le \cdots \le a_m \le b_m \le a_{m+1}
         \quad &(d = 1).
\end{aligned}
\end{equation*}
A symbol $(A,B) \in Y_{n,d}^r$ for $d \ge 0$ is said to be  
distinguished if it is an image of a distinguished symbol in 
$X_{n,d}^{r,0}$. 
\par
Assume that $r \ge 1$.  The two symbols 
$(A,B), (A', B') \in X_{n}^{r,s}$ are said to be similar if 
$A \cup B = A' \cup B', A \cap B = A' \cap B'$, namely under some
shift, $A \cup B$ coincides with $A' \cup B'$ with multiplicities.
This defines an equivalence relation on the set $X_n^{r,s}$, and 
an equivalence class is called a similarity class in $X_n^{r,s}$.
A similarity class in $Y_{n, \ \even}^{r}$ or $Y_{n, \ \odd}^r$
is defined in a similar way. 
A similarity class in $Y_{n,\ \even}^r$ 
containing $(A,A)$ is called a degenerate
class, which consists of two copies of $(A,A)$.
It is easy to see that each (non-degenerate) similarity class 
contains a unique distinguished symbol.  
\par
It is known by [L1], [LS] that 
a similarity class in $X_n^{r,s}$, $Y_{n,\ \even}^r$ or 
$Y_{n,\ \odd}^r$ is in a natural way regarded as a vector space
over $\BF_2$ as follows.  
Let $\vL = (A, B)$ be a distinguished symbol in a similarity class $\ZC$
in $X_n^{r,s}, Y_{n,\ \even}^r$ or $Y_{n,\ \odd}^r$.  We assume that 
$A \ne B$ if $(A,B) \in Y_{n,\ \even}^r$, and 
put $S = (A \cup B) \backslash  (A\cap B)$.  Then $S \ne \emptyset$
and it is written as $S = \{ c_1, c_2, \dots, c_t \}$ in an increasing 
order. A non-empty subset $I = \{ c_i, c_{i+1}, \dots, c_j\}$ of 
$S$ is called an interval if $c_{k+1} - c_k < r+s$ for $i \le k < j$
and it is maximal with respect to this condition.  
We say that $c_i$ is the tail of $I$.  An interval is called 
an initial interval if $c_i < s$.
Hence the initial interval exists only in the case where $s > 0$, 
and in that case, it exists uniquely after some shift. 
$S$ is a disjoint union of intervals.  
\par
Assume that $\ZC \subset X_{n}^{r,s}$.
Let $I$ be an interval which 
is not initial with the tail $c$.  If $c \in A$ 
(resp. $c \in B$), then there exists a unique $(A',B') \in \ZC$
such that $c \in B$ (resp. $c \in A$) and that 
$A \cap J = A' \cap J, B \cap J = B' \cap J$ for all other intervals
$J$.  This means that $(A',B')$ is obtained from $(A, B)$ by permuting 
the entries in the interval $I$.  
All the symbols in $\ZC$ are obtained from $(A,B)$ 
by permuting the entries in certain intervals. 
Let $\CI$ be the set of non-initial intervals in $S$ and $\CP(\CI)$ 
the set of all subsets of $\CI$.  The above argument shows that $\ZC$ is 
in bijection with the set $\CP(\CI)$, which has a natural structure 
of $\BF_2$ vector space with origin $\vL$  and is denoted 
by $V_{\vL}^{r,s}$.
In the case where $\ZC \subset Y_{n,\ \even}^r$ or 
$\ZC \subset Y_{n,\ \odd}^r$ ($\ZC$: non-degenerate), 
$\ZC$ is in bijection with the quotient set of $\CP(\CI)$ 
under the relation  
$\CK \sim \CI \backslash \CK$ for $\CK \in \CP(\CI)$.  Hence $\ZC$ is 
identified with the $\BF_2$ vector space $V_{\vL}^{r,0}/L$, where
$L$ is a line generated by $\CI \in \CP(\CI)$, which we denote by 
$V_{\vL}^r$.
\para{3.5.}
Let $G$ be as in 3.1. We associate the sets 
$X_n^{1,1}, Y_{n,\ \even}^2, Y_{n,\ \odd}^2$ for 
$Sp_{2n}, SO_{2n}, SO_{2n+1}$ with $p \ne 2$ and 
$X_n^{2,2}, Y_{n,\ \even}^4$ for 
$Sp_{2n}, SO_{2n}$ with $p = 2$.
Recall $n_0$ in 3.3. 
\par\medskip
(a)\ The case $X_n^{1,1}$.  We have $r = s = 1$, and
$n_0 = \frac 12 d(d-1)$ for odd $d$. Hence (3.3.3) implies 
a bijection 
\begin{equation*}
X_n^{1,1} \lra \coprod_{\substack{ d \in \BZ \\ d\ \odd}}
         (W_{n-\frac12 d(d-1)})\wg.
\end{equation*}
\par\medskip
(b)\ The case $Y_{n,\ \even}^2$ or $Y_{n,\ \odd}^2$.
We have $r = 2, s = 0$, and $n_0 = \frac 12 d^2$ for even $d$ 
and $n_0 = \frac 12 (d^2-1)$ for odd $d$.  Hence (3.3.4) implies
bijections 
\begin{align*}
Y_{n,\ \odd}^2 &\lra \coprod_{\substack{ d\ge 1 \\ d\ \odd}}
            (W_{n-\frac12 (d^2-1)})\wg, \\
Y_{n,\ \even}^2 &\lra (W_n')\wg \coprod\bigl(
       \coprod_{\substack{ d > 0 \\ d\ \even}}
                 (W_{n - \frac12 d^2)})\wg\bigr).
\end{align*}
\par\medskip
(c)\ The case $X_n^{2,2}$.  We have $r = s = 2$, and $n_0 = d(d-1)$
for odd $d$.  Hence (3.3.3) implies a bijection 
\begin{equation*}
X_n^{2,2} \lra \coprod_{\substack{ d \in \BZ \\ d\ \odd}}
         (W_{n- d(d-1)})\wg.
\end{equation*}
\par\medskip
(d)\ The case $Y_{n,\ \even}^4$.  We have $r = 4, s= 0$, and 
$n_0 = d^2$ for even $d$. 
Hence (3.3.4) implies a bijection 
\begin{equation*}
Y_{n,\ \even}^4 \lra (W_n')\wg \coprod\bigl(
       \coprod_{\substack{ d > 0 \\ d\ \even}}
                 (W_{n - d^2})\wg\bigr).
\end{equation*}
\para{3.6.}  
Let $X = X_n^{1,1}, Y_{n,\ \odd}^2, Y_{n,\ \even}^2$ or 
$X_n^{2,2}, Y_{n,\ \even}^4$ according as 
$G = Sp_{2n}$, $SO_{2n+1}$, $SO_{2n}$ with $p \ne 2$, or 
$G = Sp_{2n}, SO_{2n}$ with $p = 2$.
In view of the bijections in 3.5 and the discussion in 3.1, 
the generalized Springer correspondence can be described 
by giving a bijection between $\CN_G$ and $X$. 
By [L1], [LS], this bijection is given explicitly in such a way 
that the set $G\uni\ssim$ of unipotent classes in $G$ is 
in bijection with the
set $X\ssim$ of similarity classes in $X$.
In what follows, we define a map $\r: G\uni\ssim \to X\ssim$ by associating a
distinguished symbol $\vL = \r(C) \in X$ for each unipotent class $C$ in $G$. 
\par\medskip
(a)\ $G = Sp_{2n}$ with $p \ne 2$.  
Let $C_{\la}$ be a unipotent class of $G$ as in 2.2, where $\la$ is 
a partition of $2n$.  We express $\la$ as 
$\la_1 \le \la_2 \le \cdots \le \la_{2m}$ for some $m$, by allowing 0 
in the entries if necessary. We divide the sequence 
$\{ \la_1, \la_2, \dots, \la_{2m}\}$ into the union of blocks as 
follows.
If $\la_i$ is even, let $\{ \la_i\}$ be a block. If $\la_i = h$ is
odd, the sequence $A_h = \{ \la_k \mid \la_k = h\}$ consisting of 
even elements, 
which we write as $\{ \la_a, \la_{a+1}, \dots, \la_b \}$ for
some $b > a$.  Then we divide $A_h$ into a disjoint union of two
elements blocks 
$\{ \la_a, \la_{a+1}\}\cup \{ \la_{a+2}, \la_{a+3}\}\cup 
     \dots \cup \{ \la_{b-1}, \la_b\}$.
We define a sequence $\nu_1, \dots, \nu_{2m} $ as follows. 
Put 
\begin{equation*}
\begin{cases}
\nu_i = \la_i/2 + i &\quad\text{ if $\{ \la_i\}$ is a block}, \\ 
\nu_i  = \nu_{i+1} = (\la_i+1)/2 + i
                      &\quad\text{ if $\{ \la_i, \la_{i+1}\}$ is a block}, 
\end{cases}
\end{equation*}
and put $A = \{ 0, \nu_2, \nu_4, \dots, \nu_{2m}\}$, 
$B = \{ \nu_1, \nu_3, \dots, \nu_{2m-1}\}$.  Then $\vL = (A, B)$ gives rise
to a distinguished symbol in $X^{1,1}_n$, which is independent 
of the choice of $m$, and $C_{\la} \mapsto \vL$
gives the required bijection $\r$.
Actually, the map $\r$ was defined in [L1, 11.6].  Although the definition 
given there is not the same as ours, it is easily checked that they 
coincide with each other. 
\par\medskip
(b)\ $G = SO_N$ with $p \ne 2$. 
Let $C_{\la}$ be a unipotent class of $G$ as in 2.2, where $\la$
is a partition of $N$.  We choose $M$ such that $M \equiv N \pmod 2$
and express $\la$ as $\la_1 \le \la_2 \le \cdots \le \la_M$.
We divide the sequence $\{ \la_1, \la_2, \dots, \la_M\}$ into the 
union of blocks as follows.  If $\la_i$ is odd, let $\{ \la_i\}$
be a block.  If $\la_i = h$ is even, the sequence 
$A_h = \{ \la_k \mid \la_k = h\}$ consists of even elements.  As in the case
(a), we divide $A_h$ as a disjoint union of two elements blocks,
$\{\la_a, \la_{a+1}\} \cup\{ \la_{a+2}, \la_{a+3}\}\cup \cdots
\cup \{ \la_{b-1}, \la_b\}$.
We define a sequence $\nu_1, \nu_2, \dots, \nu_M$ as follows.
Put
\begin{equation*}
\begin{cases}
\nu_i = (\la_i-3)/2 + i &\quad\text{ if $\{ \la_i\}$ is a block}, \\
\nu_i = \nu_{i+1} = (\la_i -2)/2 + i 
            &\quad\text{ if $\{ \la_i, \la_{i+1}\}$ is a  block},
 \end{cases}
\end{equation*}
and put 
$A = \{ \nu_1, \nu_3, \dots, \nu_{[(M+1)/2]}\}$, 
$B = \{ \nu_2, \nu_4, \dots, \nu_{[M/2]}\}$.
Then $\vL = (A,B)$ gives rise to a distinguished symbol in 
$Y_{n,\ \odd}^2$ or $Y_{n,\ \even}^2$ according as $N$ is 
odd or even, which is independent of the choice of $M$. 
The map $C_{\la} \mapsto \vL$ gives the required
bijection $\r$.  The proof follows, as in the case (a), 
 from the discussion in [L1, 11.7].
\par\medskip
(c)\ $G = Sp_{2n}$ with $p = 2$.
The map $\r$ is defined in [LS, 2.1]. The following definition 
is slightly modified from the original one so as to 
fit to the case (a). 
Let $C_{\la,\ve}$ be a unipotent class of $G$ as in 2.3, where
$\la$ is a partition of $2n$.
Here $\la$ is the same as in the case (a), and we 
express it as $\la_1 \le \cdots \le \la_{2m}$ for some $m.$
We use the convention that $\ve(0) = 1$.
We divide the 
set $\{\la_1, \la_2,  \dots, \la_{2m}\}$ into a disjoint
union of blocks as follows.  If $\ve(\la_i) = 1$, then 
$\{ \la_i\}$ is a block.  If $\ve(\la_i) = 0$ or $\w$ for 
$\la_i = h$, the sequence $A_h = \{ \la_k \mid \la_k = h\}$
has even cardinality, and it is divided into blocks 
as in the case (a).
We define a sequence $\nu_1, \dots, \nu_{2m}$ as follows. 
\par
(i)\ If $\{\la_i\}$ is a block, put
\begin{equation*}
\nu_i = \la_i/2 + 2i.
\end{equation*}
\par
(ii)\ If $\{\la_i, \la_{i+1}\}$ is a block and $\ve(\la_i) = \w$,
put
\begin{align*}
&\nu_i = (\la_i + 1)/2 + 2i, \\
&\nu_{i+1} = \nu_i + 1.  
\end{align*}
\par
(iii)\ If $\{ \la_i, \la_{i+1}\}$ is a block and $\ve(\la_i) = 0$, 
put
\begin{align*}
&\nu_i = (\la_i+2)/2 + 2i, \\
&\nu_{i+1} = \nu_i.
\end{align*} 
Put $A = \{ 0, \nu_2, \nu_4, \dots, \nu_{2m}\}$,
$B = \{ \nu_1, \nu_3, \dots, \nu_{2m-1}\}$.  Then by [LS, 2.2], 
$\vL = (A, B)$ gives rise
to a distinguished symbol in $X_n^{2,2}$, which is independent of the
choice $m$, and $C_{\la,\ve} \mapsto \vL$ gives the required bijection 
\par\medskip
(d)\ $G = SO_{2n}$ with $p = 2$.
The map $\r$ is defined in [LS, 3.1].  
Let $C_{\la,\ve}$ be a unipotent class in $G$.
$\la$ is the same as in the case (c), and we define the 
block in the same way as the case (c).  However, here 
we use the convention that $\ve(0) = 0$.  Note that in the sequence
$\{ \la_1, \dots, \la_{2m}\}$ the multiplicity 
of 0 is even since the number of non-zero $\la_i$ is even (cf. 2.3).
We define a sequence $\nu_1, \dots, \nu_{2m}$ as follows.
\par
(i)\ If $\{\la_i\}$ is a block, put
\begin{equation*}
\nu_i = (\la_i-6)/2 + 2i.
\end{equation*}
\par
(ii)\ If $\{ \la_i, \la_{i+1}\}$ is a block and $\ve(\la_i) = \w$, 
put
\begin{align*}
&\nu_i = (\la_i-5)/2 + 2i, \\
&\nu_{i+1} = \nu_{i} +1.
\end{align*}
\par
(iii)\ If $\{ \la_i, \la_{i+1}\}$ is a block and $\ve(\la_i) = 0$, 
put
\begin{align*}
&\nu_i = (\la_i - 4)/2 + 2i, \\
&\nu_{i+1} = \nu_i. 
\end{align*}
Put $A = \{ \nu_1, \nu_3, \dots, \nu_{2m-1}\}$,
$B = \{ \nu_2, \nu_4, \dots, \nu_{2m}\}$.  Then 
by [LS, 3.2],  
$\vL = (A,B)$ gives rise to a distinguished symbol in 
$Y_{n,\ \even}^4$, which is independent of the choice of $m$, 
and $C_{\la,\ve} \mapsto \vL$ gives the required bijection. 
\para{3.7.}
We return to the setting in the beginning of 3.6, and 
let $\r : G\uni\ssim \to X\ssim$ be the bijection constructed 
in 3.6.  By making use of $\r$, we shall construct a bijection
$\wt\r : \CN_G \to X$.  For a unipotent class $C$ in $G$, 
let $\ZC$ be the similarity
class in $X$ containing the distinguished symbol $\vL = \r(C)$.
Take $u \in C$. As discussed in 3.4, $\ZC$ has a natural structure of
$\BF_2$-vector space $V_{\vL}^{r,s}$ for $X = X_n^{1,1}$ or $X_n^{2,2}$
with the basis corresponding to the set of intervals.  Since 
$A_G(u)$ is an elementary abelian 2-group, it has a natural structure of
$\BF_2$-vector space, and so does the dual group 
$A_G(u)\wg$. It was 
shown in [L1, 11], [LS, 2.2] that $V_{\vL}^{r,s}$ is naturally
identified with $A_G(u)\wg$, where the set of intervals is in bijection
with the set of generators in $A_G(u)$ given in 2.4;  if $I$ is an 
interval corresponding to the generator $a_i$ of $A_G(u)$, we associate 
the character $\x_i$ of $A_G(u)$ which takes the value $-1$ for $a_i$ and 
1 for other generators.  Similar argument also works for the case where
$X = Y_{n,\ \odd}^2, Y_{n,\ \even}^2$ and $Y_{n,\ \even}^4$, and by 
[L1, 11], [LS, 3.2] $V_{\vL}^r$ is naturally identified with
$A_G(u)\wg$. 
Hence the similarity class $\ZC$ can be identified with $A_G(u)\wg$.
\par
On the other hand, for a fixed $C$, $G$-equivariant simple
local systems on $C$ are parametrized by $A_G(u)\wg$.
Thus combining with the above argument, we obtain a bijection 
$\wt\r : \CN_G \to X$.  This bijection describes combinatorially 
the generalized Springer correspondence, namely, 
\begin{thm}[{[L1, 12.3,13.3], [LS, 2.4, 3.3]}]
Let $G$ be as in 3.1.  Then the composite of $\wt\r$ with the 
bijection in 3.5 gives the generalized Springer correspondence 
in 3.1.
\end{thm} 
\section{Mail results}
\para{4.1.}
Let $G$ be as in 2.1 and we apply the argument in Section 1 for $G$.
First consider the case where $G$ is of split type.  For each 
unipotent class 
$C'$ in $G$, we choose a split element $u \in {C'}^F$ described in 
Section 2, i.e., (2.9.3) for $G = Sp_{2n}$ or $SO_N$ with $p \ne 2$, 
(2.7.3) for $G = Sp_{2n}$ with $p = 2$, (2.8.1) for $G = SO_N$ with
$p = 2$. For each pair $(C',\CE') \in \CN_G^F$, we fix 
an isomorphism $\p_0 = \p_{\CE'}: F^*\CE' \to \CE'$ in 1.6 so that 
the induced isomorphism 
$\CE'_u \to \CE'_u$ is identity. (Note that $F$ acts trivially on 
$A_G(u)$.  This is known by [SS] for the case where $p \ne 2$, and 
follows from (2.7.4) for the case where $p = 2$.) 
Take $(L, C,\CE) \in \CM^F_G$.  Then $L$ is of the same 
type as $G$ of split type.  Thus we choose the split 
element $u_0 \in C^F$.  We fix an isomorphism 
$\vf_0: F^*\CE \to \CE$ in 1.6 so that the induced isomorphism
$\CE_{u_0} \to \CE_{u_0}$ is identity.  Let $V_{(C',\CE')}$ be the 
irreducible $\CW$-module and $\s_{(C',\CE')}$ be the isomorphism 
on $V_{(C',\CE')}$ given in 1.6.  Since $F$ is of split type, 
$F$ acts trivially on $\CW$ and so $\s_{(C', \CE')}$ commutes with 
the action of $\CW$.  It follows that $\s_{(C',\CE')}$ is a scalar
map.  
\par
Next consider the case where $G = SO_{2n}$ with $F = F_0\s$ 
of non-split type.
We choose the split elements $u' \in {C'}^F$ and  $u_0' \in C^F$ as 
in 2.10, and fix $\p_0 = \p_{\CE'}, \vf_0$ as above. 
(Again $F$ acts trivially on $A_G(u')$ by (2.10.2).)
Let 
$(L,C,\CE) \in \CM_G^F$. If $L \ne T$, 
$\CW$ is a Weyl group of type $B$ and so $F$ acts trivially on $\CW$.  
Hence $\s_{(C',\CE')}$ is a scalar map.
While if $L = T$, $\CW = W_n'$ is the Weyl group of type $D_n$ 
and $F$ acts non-trivially on $\CW$. Note that $\s$ acts on 
$W_n'$ and the semidirect product $W_n'\rtimes\lp\s\rp$ is isomorphic to 
$W_n$.  Assume that $V_{(C',\CE')} = E$ is $F$-stable.  Then
$E$ can be extended to an irreducible representation $\wt E$ of $W_n$
via the map $\s_{(C',\CE')}$. Since $\s\iv \circ \s_{(C',\CE')}$ 
commutes with the action of $\CW$, it acts as a scalar on $E$. 
Thus in order to describe the map $\s_{(C',\CE')}$, we have only to
determine this scalar together with the representation $\wt E$.
\par
We recall the preferred extension $\wt E$ of $W_n$ due to 
[L2, IV, 17.2].
An $F$-stable irreducible representation $E$ of $W_n'$ is parametrized 
by an unordered pair $(\a; \b)$ of partitions such that 
$|\a| + |\b| = n$ and that $\a \ne \b$.  
We write $\a : \a_1 \le \cdots \le \a_m$, 
$\b: \b_1 \le \cdots \le \b_m$, and define a symbol 
$\vL_E$ associated to $E$ by an unordered pair $(\la; \mu)$ 
with $\la_i = \a_i + i-1, \mu_j = \b_j + j-1$. (This is a different
type of symbols from those appeared in Section 3.) Irreducible 
representations of $W_n$ are parametrized in a similar way, but
by using an ordered pair $(\a; \b)$ and its associated symbol
$(\la;\mu)$.  For an $F$-stable irreducible representation $E$, 
there exists two extensions to $W_n$, which correspond to two
symbols $(\la; \mu)$ and $(\mu; \la)$ for $W_n$.  
An extension $\wt E$ of $E$ is called the preferred extension of
$E$  if in the symbol $\vL_{\wt E}$, the smallest number which 
does not appear  
in both entries appears in the second entry.  For example, 
$(n; 0)$ is the symbol associated to the unit representation 
of $W_n'$, and it is extended 
to the unit representation $(n; 0)$ or the long sigh representation 
$(0; n)$ of $W_n$. 
In this case, $(n; 0)$ is the preferred extension. 
\par
We can state our main results.
\begin{thm}  
Let $G = Sp_{2n}, SO_N$ with $F$ of split type ($N$ is even if $p = 2$).
Then $\s_{(C',\CE')}$ is $q^{(a_0+r)/2}$ times identity.
\end{thm}
\begin{thm}   
Let $G = SO_{2n}$ with $F$ of non-split type. 
\begin{enumerate}
\item
Suppose $L \ne T$. Then 
$\s_{(C',\CE')}$ is $q^{(a_0 + r)/2}$ times identity.
\item 
Suppose $L = T$.  Then $\s_{C', \CE'} = q^{(a_0 + r)/2} \s$, 
and $W'_n\lp\s\rp$-module $\wt E$ coincides with the preferred extension
of $E$. 
\end{enumerate}
\end{thm}
In view of Lemma 1.11, we have the following corollary.
\begin{cor} 
Let $G$ be as in 2.1.  For each $(C',\CE') \in \CN_G^F$, 
choose $\p_0: F^*\CE' \isom \CE'$ by choosing a split element 
$u \in {C'}^F$.  Let $\g$ be the constant as in 1.10.  
Then we have $\g = 1$, namely we have $Y_j = Y_j^0$ for 
$j = (C',\CE')$.   
\end{cor}
 \para{4.5.}  We prove the theorem by making use of the restriction 
formula in Corollary 1.9.  We choose a standard parabolic subgroup 
$Q = MU_Q$ 
such that the Levi subgroup $M$ of $Q$ is of the same type as $G$ 
with semisimple rank $n-1$.
Take $(L,C,\CE) \in \CM_G^F$ and let $P$ be the $F$-stable standard 
parabolic subgroup of $G$ whose Levi subgroup is $L$.  Then we have
$P \subset Q$ and $L \subset M$.
Take $(C',\CE') \in \CN_G^F$, and take a split element 
$u \in {C'}^F$.  Also take $(C_1, \CE_1) \in \CN_M^F$, and choose
a split element $v \in C_1^F$.  Let $\r$ (resp. $\r_1$) be the irreducible
representation of $A_G(u)$ (resp. $A_M(v)$) corresponding to $\CE'$
(resp. $\CE_1$). Since $F$ acts trivially on $A_G(u)$
and $A_M(v)$, 
the extension $\wt{ \r\otimes \r_1^*}$ to 
$\wt A(u,v)$ in 1.8 is just the trivial extension of $\r\otimes \r_1^*$.
Then we have the following lemma.
\begin{lem} 
Assume that $F$ acts trivially on $\CW$ and on $\CW_1$.
Let $E \in \CW\wg$ (resp. $E_1 \in \CW_1\wg$) be corresponding to 
$(C',\CE') \in \CN_G$ (resp. $(C_1, \CE_1) \in \CN_M$), and assume that   
$E_1$ occurs in the restriction of $E$ to $\CW_1$.  Suppose that the
theorem holds for $\s_{(C_1, \CE_1)}$.  If $X_{u,v} \ne \emptyset$ and 
$F$ acts trivially on $X_{u,v}$, then the theorem holds for 
$\s_{(C',\CE')}$.
\end{lem}
\begin{proof}
We follow the notation in Section 1.   
Since $E_1$ occurs in the restriction of $E$ to $\CW_1$, 
$M_{E_1} \ne 0$ in (1.6.1). 
Since $\s_{(C',\CE')}$ and $\s_{(C_1,\CE_1)}$ is a scalar map, 
$\s_{\CE_1,\CE'}$ is also a non-zero scalar map by (1.6.2).  
Since $F$ acts trivially on $X_{u,v}$, $\wt\ve_{u,v}$ is the
trivial extension of 
$\ve_{u,v}$ to $\wt A(u,v) = \lp\t\rp\times  A(u,v)$.
It follows that 
\begin{equation*}
\lp\wt\ve_{u,v}, \wt{\r\otimes \r_1^*}\rp_{A(u,v)\t}
    = \lp\ve_{u,v}, \r\otimes \r_1^*\rp_{A(u,v)}.
\end{equation*}
Then Corollary 1.9 together with Corollary 1.5 implies that 
\begin{equation*}
\Tr(\s_{\CE_1,\CE'}, M_{E_1}) = 
          q^{-d_{C_1,C'} + \dim U_Q}\dim M_{E_1}, 
\end{equation*}
and we see that $\s_{\CE_1,\CE'}$ is a scalar map 
by $q^{-d_{C_1,C'} + \dim U_Q}$. 
By our assumption, $\s_{(C_1,\CE_1)}$ is a scalar map by 
$q^{(a_0' + r')/2}$, where $a_0', r'$ are as in 1.10 with respect to
$M$.
Thus again by (1.6.2), we see that $\s_{(C',\CE')}$ is a scalar map by 
$q^{(a_0 + r)/2}$. 
\end{proof}
\para{4.7.} 
In view of Lemma 4.6, it is important to know the Frobenius action 
on $X_{u,v}$.
We note that
\par\medskip\noindent 
(4.7.1)\  $Z_G(u) \times Z_M(v)U_Q$ acts transitively on 
the set $Y_{u,v}$.  
\par\medskip
In fact, put 
$\CQ_{u, C_1} = \{ gQ \in G/Q \mid g\iv ug \in C_1U_Q\}$.  $\CQ_{u,C_1}$
is a locally closed subvariety of $G/Q$.
We have a surjective morphism $Y_{u,v} \to \CQ_{u,C_1}, g \mapsto gQ$, 
which induces an isomorphism $Y_{u,v}/Z_M(v)U_Q \simeq \CQ_{u,C_1}$.
It is known by [Sp1, II, 6.7] that $Z_G(u)$ acts
transitively on $\CQ_{u,C_1}$.
(4.7.1) follows from this. 
\begin{lem} 
If $Y_{u,v}^F \ne \emptyset$, then $F$ acts trivially on 
$X_{u,v}$. 
\end{lem}
\begin{proof}
\par
By (4.7.1), the closure of $Z_G^0(u)gZ_M^0(v)U_Q$ in $Y_{u,v}$ 
($g \in Y_{u,v}$) gives an element $x \in X_{u,v}$, and 
$A_G(u)xA_M(v)$ gives all the irreducible components in $Y_{u,v}$.
We can choose $g \in Y_{u,v}^F$.  Hence 
$Z_G^0(u)gZ^0_M(v)U_Q$ is $F$-stable, and so is $x$.  
Since $F$ acts trivially on $A_G(u)$ and $A_M(v)$, 
$F$ stabilizes all the irreducible components in $Y_{u,v}$.  
\end{proof}
\para{4.9.}  
Assume that $F$ is of split type. 
For each $u \in G^F$, we shall give $v \in M^F$ such that 
$Y_{u,v}^F \ne \emptyset$.  Let $V$ and $f$ be as in 2.1.
Note that $G/Q$ can be identified with the subvariety of $\PP(V)$
consisting of $\lp x\rp$ for isotropic vectors $x$ with respect to $f$ 
($\lp x\rp$ denotes the line spanned by $x$).
Under the setting in 2.7, 2.9, we consider the vector space 
$\bar V = \bigoplus M_j$ which is identified with $V^F$.
In the following cases, we can find $g \in G^F$ such that 
$g\iv ug \in Q$ and that 
$v = \pi(g\iv ug)$ is a split element in $M$  
($\pi: Q \to M$ is the natural projection). 
In particular, we have $g \in Y_{u,v}^F$. 
In the discussion below, we identify the partitions with 
the corresponding Young diagrams.
\par  First we consider the case where $G = Sp_{2n}$ or 
$SO_N$ with $p \ne 2$.
\par\medskip
(i)\ Take $M_j$ such that $\ve(\la_j) = 1$. 
Let $e_1^j, \dots, e_{h}^j$  be the basis of $M_j$ 
with $h = \la_j$ given in 2.9 (a).  The stabilizer 
of $\lp e_1^j\rp$ in $G$ is
a parabolic subgroup $gQg\iv$ with some $g \in G^F$. 
The nilpotent transformation $X_j$ on $M_j$ induces a
map $\bar X_j$ on $\bar M_j = \lp e_1^j\rp\orth/\lp e_1^j\rp$,
and one can define a nilpotent element $\bar X$ on 
$\bar M_j \oplus \bigoplus_{j' \ne j}M_{j'}$.  This determines 
a unipotent element $\pi(g\iv ug) = v$ in $M^F$ 
which is a split element.
In this case $v$ is of type $\la'$, where $\la'$ is obtained from 
$\la$ by replacing one row $h$ such that $\ve(h) = 1$ by $h - 2$.
\par
(ii)\ Take $M_j$ such that $\ve(\la_j) = 1$ and that  
$c_h$ is even for $h = \la_j$. We choose $j$ such that 
$\la_j = \la_{j-1} = h$ and 
consider $N_j = M_j \oplus M_{j-1}$ with $f_j' = f_j + f_{j-1}$.
 Thus $N_j$ has a basis 
$e_1^j, \dots, e_h^j, e_1^{j-1}, \dots, e_h^{j-1}$.
By our construction, we have 
\begin{equation*}
f_j(e_1^j, e_h^j) = - f_{j-1}(e_1^{j-1}, e_h^{j-1}) = \pm 1.
\end{equation*}
We consider 
$\bar N_j = \lp x\rp\orth/\lp x\rp$ for $x = e_1^j + e_1^{j-1}$.
The nilpotent transformation $X_j + X_{j-1}$ on $M_j \oplus M_{j-1}$
induces a linear map $\bar X_j$ on $\bar N_j$, and one can define 
a nilpotent element $\bar X$ on $\bar N_j \oplus\bigoplus M_{j'}$.
This determines a unipotent element $\pi(g\iv ug) = v$ in $M$ of type
$\la'$, where $\la'$ is obtained from $\la$ by replacing two rows of 
length $h$ by two rows of length $h-1$.  It is easy to see that 
$v$ is a split element in  $M^F$.
\par
(iii)\ Take $M_j$ such that $\ve(\la_j) = \w$.  Put $h = \la_j$. 
We choose a basis of $M_j$ as in 2.9 (b).  
Let  $x = e_1^j + e_1^{j-1}$. 
Then $x$ is an isotropic vector with respect to $f_j$.
Now $X_j$ induces $\bar X_j$ on 
$\bar M_j = \lp x\rp\orth/\lp x \rp$, and one can define 
$\bar X$ on $\bar M_j \oplus \bigoplus_{j' \ne j}M_{j'}$.
This determines a unipotent element 
$v = \pi(g\iv ug) \in M^F$.  In this case, $v$ is of type 
$\la'$, where $\la'$ is obtained from $\la$ by replacing 
two rows of length $h$ by two rows of length $h - 1$.  
It is easy to see that $v$ is a split element in  $M^F$.
\par\medskip
It is known (cf. [S1]) that if $Y_{u,v} \ne \emptyset$, 
for $u \in C_{\la}, v \in C_{\la'}$, then $\la'$ is obtained
from $\la$ by the procedure described above.  Hence if $v$ is 
split, it coincides with one of the above cases.  It follows that 
\begin{prop} 
Let $G = Sp_{2n}$ or $SO_N$ with $p \ne 2$.  Assume that 
$F$ is of split type.  Let $u \in G^F$ 
and $v \in M^F$ be such that $Y_{u,v} \ne \emptyset$.
Assume that $u,v$ are split elements.  Then $F$ acts trivially on 
$X_{u,v}$.
\end{prop}
\para{4.11.}
Next we consider the case where 
$G =  Sp_{2n}$ with $p = 2$.  We keep the setting in 4.9, 
in particular assume that $F$ is of split type.
\par
(i)\ Choose $M_j$ as in 2.7 (a), i.e., the case where 
$\ve(h) = 1$ for $h = \la_j$. 
$v_j$ induces a linear map $\bar v_j$ on 
$\bar M_j = \lp e_1^j\rp\orth/\lp e_1^j\rp$ and one can define
$\bar v$ on $\bar M_j \oplus\bigoplus_{j' \ne j}M_{j'}$.  
This determines a unipotent element $\pi(g\iv ug) = v \in M$ for some 
$g \in G^F$.  
$v$ is contained in the class 
$C_{\la',\ve'}$ in $M$, where $\la'$ is obtained from $\la$ by replacing
one row of length $h$ by a row of length $h-2$, and $\ve'$ is given by 
$\ve'(h-2) = 1$ and $\ve'(\la_k') = \ve(\la_k)$ if $\la_k' \ne h-2$.
By our construction of the form $f_j$ in 2.5, 
we see that $v$ is a split element if $h-2$ does not occur in the row of 
$\la$, nor $h-2$ occurs and $\ve(h-2) = 1$.  
\par
(ii)\ Assume that $c_h \ge 2$ and $\ve(h) = 1$ 
for $h = \la_j$, and take 
$\la_j = \la_{j-1} = h$.  We consider $M_j = V_j \oplus V_{j-1}$ 
with $f_j = f_j^0 + f_{j-1}^0$, 
where $f_j^0, V_j$ are as in 2.7. Let $x = e_1^j + e_1^{j-1}$.  
Put $\bar M_j = \lp x\rp\orth/\lp x \rp$.  Then $v_j + v_{j-1}$
induces a unipotent element $\bar v_j$ on $\bar M_j$.
This determines $\pi(g\iv ug) = v$ of $M^F$ as before.  This 
construction is exactly the same as the one in 2.7 (b).
Hence $v$ is  contained in $C_{\la',\ve'}$ 
in $M$, where 
$\la'$ is obtained from $\la$ by removing two rows of length $h$
by two rows of length $h-1$, and $\ve'(h-1) = \w$, 
$\ve'(h') = \ve(h')$ for all $h' \ne h$.
In this case, $v$ is a split element without any condition. 
\par
(iii)\ Choose $M_j$ as in 2.7 (b), i.e., the case where $\ve(h) = \w$.
The basis of $M_j$ is given in 2.7 (b).  Then 
$\bar M_j = \lp \bar e_1^j\rp\orth/\lp \bar e_1^j\rp$ has a basis 
$\bar e_{n-1}^j, \bar e_{n-1}^{j-1}, \dots, \bar e_2^j, \bar e_2^{j-1}$
(we use the same notation for the image on $\bar M_j$ as the one in
$M_j$).
Thus the induced linear map $\bar v_j$ on $\bar M_j$ is just a sum of 
two copies of nilpotent elements as given in 2.7 (a) with respect to the 
induced form $\bar f_j$ on $\bar M_j$.
It follows that $\bar v$ determines a unipotent element 
$\pi(g\iv ug) = v \in M$.  $v$ is contained in $C^F_{\la', \ve'}$, where 
$\la'$ is obtained from $\la$ by replacing two rows of length $h$ by 
two rows of length $h-1$. $\ve'$ is given by $\ve'(h-1) = 1$,
and is the same as $\ve$ for all other $h'\ne h$.  
In particular, $v$ is a split element if $h-1$ does not occur in the rows
of $\la$ nor $\ve(h-1) = 1$.
\par
(iv)\ Choose $M_j$ as in (iii).  Under the notation there, put 
$x = \bar e_2^j + \bar e_2^{j-1}$. 
Then $\bar M_j = \lp x\rp\orth/\lp x \rp$ has a basis
exactly the same as the basis of $N\orth/N$ in 2.7 (c).  Thus
the induced map $\bar v_j$ on $\bar M_j$ determines a unipotent 
element $\pi(g\iv ug) = v$ of $M$ in the same way as above.
$v$ is a split element in $M$ and is contained in $C_{\la',\ve'}$, 
where $\la'$ is obtained from $\la$ by replacing two rows of length 
$h$ by two rows of length $h-1$, and $\ve'(h-1) = 0$ if $h-1$ does not
occur in $\la$ nor if $\ve(h-1) = 0$.
\par
(v)\ Choose $M_j$ as in 2.7 (c), i.e., the case where $\ve(h) = 0$.
The basis of $M_j$ is given in 2.7 (c).  Then 
$\bar M_j = \lp \bar e_1^j\rp\orth/\lp\bar e_1^j\rp$ has a basis
$\bar e_{n-1}^j+ \bar e_{n-1}^{j-1}, 
    \bar e_{n-2}^j, \bar e_{n-2}^{j-1}, \dots, 
      \bar e_3^j, \bar e_3^{j-1}$,
     $\bar e_2^j = \bar e_2^{j-1}$.
Thus the induced linear map $\bar v_j$ on $\bar M_j$ is the same as 
the case (b) as above with respect to the induced form $\bar f_j$.
It follows that $\bar v$ determines a unipotent element 
$\pi(g\iv ug) = v \in M$. $v$ is contained in $C_{\la',\ve'}$
where $\la'$ is obtained from $\la$ by replacing two rows of length
$h$ by two rows of length $h-1$, and $\ve'$ is given by 
$\ve'(h-1) = \w$ if $h-1$ does not occur in the rows of $\la$.
In this case, $v$ is a split element without any 
 condition.  
\par\medskip
Finally, we consider the case where $G = SO_{2n}$ with $p = 2$.
The argument in the case of $G = Sp_{2n}$ with $p = 2$
works well for this case under a suitable modification, 
since in each case of (II), the induced quadratic 
form is of the same type
as the original one, and so one can check easily 
that $v$ is a split element.
\para{4.12.}
We consider the counter part of Proposition 4.10 and 4.11 for 
the case where $F = F_0\s$ is of non-split type.  We assume that
$G = SO_{2n}$ ($p:$ arbitrary) and $\wt G = O_{2n}$.  
Let $u \in C^{F_0}_{\la, \ve}$ and 
$v \in C^{F_0}_{\la',\ve'}$ be split elements as in 4.9 or 4.11, and 
assume that $Y_{u,v}^{F_0} \ne \emptyset$.  Take 
$g \in Y_{u,v}^{F_0}$.  Then $g\iv ug \in vU_Q$.  By replacing 
$u$ by $g\iv ug $ (an element in the split class)
we may assume that $u \in vU_Q$.
Let $a = a_i \in Z_{\wt G}(u)$ be as in (2.10.1) 
and  $\da \in Z_{\wt G}(u)$ be its representative.
Let $a' \in A_{\wt M}(v)$ be defined similar to $a$.  
We assume that $a$ and $a'$ are both related to the row
of the same length in $\la$ and $\la'$.
By the explicit description of the element $\da \in Z_{\wt G}(u)$
(see 2.7 $\sim$ 2.10), we see that $\da$ normalizes $Q$.  It follows
that $\da \in Z_{\wt M}(v)$, which gives a representative of $a'$. 
Let us write $\da = xs$ and let $\a\iv F(\a) = x$ be as in (2.10.1).
Since we may take $s \in M$, we have $x \in M$, and so 
$\a \in M$ also.
Let $u' = \a u\a\iv \in C_{\la,\ve}^F$
$v' = \a v\a\iv \in C_{\la', \ve'}^F$ be split elements as given
in (2.10.1).  We have the following lemma.
\begin{lem}  
Let the notations be as above, and $u,v$ (resp. $u',v'$) 
be split elements with respect to $F_0$ (resp. $F$).  
Then under the assumption in Proposition 4.10 or 4.11, 
$F$ acts trivially on the set $X_{u',v'}$. 
\end{lem}
\begin{proof}
Since $\a \in M$, $\ad \a$ maps $Y_{u,v}$ onto $Y_{u',v'}$, and so
induces a bijection between $X_{u,v}$ and $X_{u',v'}$.  Since 
$\da F_0 = xF$, $\ad \a$ maps $\da F_0$-stable elements of $X_{u,v}$ 
to $F$-stable elements of $X_{u',v'}$.  In view of Proposition 4.10 
and 4.11, 
we may assume that $F_0$ acts trivially on the set $X_{u,v}$.
Hence in order to prove the lemma, it is enough to show that 
any element in $X_{u,v}$ is stable by $\ad \da$. 
Now an irreducible component of $Y_{u,v}$ is expressed as the closure
of $Z_G^0(u)gZ_M^0(v)U_Q$ for some $g \in Y_{u,v}$.  By our choice 
of $u,v$, we may take $g = 1$. Then $Z_G^0(u)Z_M(v)^0U_Q$ is 
stable by $\ad \da$, and so there exists an irreducible component stable
by $\ad \da$.
Since $A_G(u) \times A_M(v)$ acts transitively 
on the set $X_{u,v}$, and $a$ commutes with $A_G(u)$ and $A_M(v)$, 
we conclude that $\ad \da$ stabilizes each element in $X_{u,v}$.
This proves the lemma. 
\end{proof}
\para{4.14.}
In order to apply Lemma 4.6, we need first to know the condition for 
$X_{u,v} \ne \emptyset$.  By the isomorphism 
$Y_{u,v}/Z_M(v)U_Q \simeq Q_{u,C_1}$ in 4.7, the elements in 
$X_{u,v}$ corresponds to the irreducible components of $Q_{u,C_1}$
of dimension $(\dim Z_G(u) - \dim Z_M(v))/2$. 
The condition for $C_1$ for the existence of such an irreducible 
component in $Q_{u,C_1}$ is described in [Sp1, II, 6.7].  By making use
of 3.6, it is interpreted in terms of the symbols 
(cf. [LS, 2.6]); 
\par\medskip\noindent
(4.14.1)\  Let 
$\vL = \r_G(u), \vL' = \r_M(v)$ be the distinguished symbols associated
to $u$ and $v$.  Then $X_{u,v} \ne \emptyset$ if and only if $\vL'$ is 
obtained from $\vL$ by decreasing one of the entries of $\vL$ by 1.    
\par\medskip
Let $E \in \CW\wg, E_1 \in \CW_1\wg$ be as in Lemma 4.6.  
In applying the lemma, we also need to know when $E_1$ appears
in the restriction of $E$.  This is given as follows (cf. [LS, 2.8]).
\par\medskip\noindent
(4.14.2)\ Let $X$ be the set of symbols as in 3.6.  Let $(A,B) \in X$
corresponding to $E \in \CW\wg$.  Then $E_1 \in \CW_1\wg$ appears 
in the restriction of $E$ if and only if the symbol 
$(A',B')$ corresponding to $E_1$ is obtained from $(A, B)$ by 
decreasing one of the entries in $A$ or $B$ by 1.
(This holds also for the case of degenerate symbols. If 
$E, E' \in \CW\wg$ are corresponding to the degenerate symbol 
$(A,A)$ and its copy, $E$ and $E'$ have the same restriction on 
$\CW_1$, and its components are parametrized by symbols obtained
by decreasing one of the entries in $(A,A)$ by 1.)
In particular, if $E_1$ appears in the restriction of $E$, we have
$X_{u,v} \ne \emptyset$.
\para{4.15.}
We are now ready to prove Theorem 4.2.  So, assume that $F$ is 
of split type.  First consider the case 
where $G = Sp_{2n}$ or $SO_N$ with $p \ne 2$.  By induction on the 
rank of $G$, we may assume that the theorem holds for 
$Sp_{2n-2}$ or $SO_{N-2}$.  Let $(C',\CE') \in \CN_G$ be corresponding 
to $E \in \CW\wg$. Let $E_1 \in \CW_1\wg$ be an irreducible component of  
the restriction of $E$, and $(C_1, \CE_1) \in \CN_M$ be the
corresponding element.  We choose the split elements 
$u \in {C'}^F$ and $v \in C_1^F$.  Then  $X_{u,v} \ne \emptyset$ 
by (4.14.2), and 
$F$ acts trivially on $X_{u,v}$ by Proposition 4.10.  Thus the assertion  
holds for $\s_{(C',\CE')}$ by Lemma 4.6, and the theorem follows.
\par
Next we consider the case where $G = Sp_{2n}$ or $SO_{2n}$ with 
$p = 2$.  Since the result in 4.11 is somewhat weaker than 
Proposition 4.10, we need a more precise argument.
Let $\la_1 \le \la_2 \le \cdots \le \la_{2m}$ be the sequence 
of $\la$ as in 3.6, and $\vL$ be the distinguished symbol associated 
to $\la$. For a given integer $h$, let 
$A_h = \{ \la_j, \la_{j+1}, \dots, \la_k\}$ be the subsequence 
consisting of $\la_i = h$.  We say that $I = \{ a, \dots, b\}$ 
($a \le b$) is a semi-interval
if $I$ corresponds to the sequence $A_h$ under the construction
of $\vL$ in 3.6 (c), (d). The element  $a$ is called the tail of 
the semi-interval $I$.
They have the following forms.  
We denote by $I_h^{\ve(h)}$ the semi-interval corresponding to $h$
and $\ve(h)$.
\begin{equation*}
\tag{4.15.1}
I_h^{\ve(h)} = \begin{cases}
              \{ a, a+2, a+4, a+6, \dots \} 
                   &\quad\text{ if } \ve(h) = 1, \\
              \{ a, a, a+4, a+4, a+8, a+8, \cdots\}
                   &\quad\text{ if }  \ve(h) = 0, \\
             \{ a, a+1, a+4, a+5, a+8, a+9, \cdots \}
                   &\quad\text{ if } \ve(h) = \w.
               \end{cases}
\end{equation*}
For two semi-intervals 
$I = \{ a, \dots, b\}, I' = \{ a', \dots, b'\}$
with $b'< a$, the distance of $I,I'$ is defined by 
$a - b'$.  It is easy to see that the distance of two 
semi-intervals is always $\ge 3$.  The case where 
the distance is 3 occurs in the following three cases;
\begin{equation*}
(I_h^1, I_{h+2}^1), \quad 
          (I_h^1, I_{h+1}^{\w}), \quad (I_h^{\w}, I_{h+1}^1).
\end{equation*}
The semi-interval $I$ is a part of an interval unless $I = I_h^0$, 
and if the distance of $I$ and $I'$ is equal to 3, they are joined to be 
a part of one big interval.  This explains the condition of generators
of $A_G(u)$ in 2.4.
\par
By induction, we may assume that the theorem holds for the smaller rank
case.  Let $E \in \CW\wg$ corresponding to $(C',\CE') \in \CN_G$.
Let $\vL$ be the distinguished symbol associated to $C'$.  We write
$C' = C_{\la,\ve}$.  Suppose that there exist two semi-intervals 
$I = \{ a, \dots, b\}, I' = \{ a', \dots, b'\}$ 
such that the distance $a - b' \ge 5$.  Then $a$ is also a tail of 
an interval, and it is easy to check that 
one can decrease $a$ by $a-1$ in $\vL$ to obtain a new symbol 
$\vL'$.  This procedure is also valid for a symbol similar to $\vL$.  
Moreover
if $a$ corresponds to $h$, and $a-1$ corresponds to $h' = h-1$ or 
$h-2$ under $C_{\la,\ve} \lra \vL$, then we have $c_{h'} = 0$.   
Let $C_{\la',\ve'} \lra \vL'$ with $v \in C_{\la',\ve'}^F$. 
There exists  
$(C_1, \CE_1) \in \CN_M$ ($C_1 = C_{\la',\ve'}$) 
corresponding to $E_1 \in \CW_1\wg$ such that 
$E_1$ occurs in the restriction of $E$.  
By making use of 4.11, we see that $Y_{u,v}^F \ne \emptyset$ for 
a split element $v \in C_{\la',\ve'}$, 
and so $F$ acts trivially on $X_{u,v}$ by Lemma 4.8. 
Now Lemma 4.6 can be applied to show that the theorem holds 
for $\s_{(C',\CE')}$.
\par
Thus it is enough to consider the case where the distance of 
$I,I'$ is $\le 4$.  Let $I = I_h^{\w}$.  There are three possibilities
for $I'$ with distance $\le 4$, i.e., $I' = I_{h-2}^{\w}, I_{h-1}^1$ or 
$I_{h-1}^0$.  For each case, one can find $(C_1, \CE_1) \lra E_1$ with 
a split element $v \in C^F_1 = C^F_{\la',\ve'}$, to which Lemma 4.6 can 
be applied. For example consider the case where $I' = I_{h-1}^0$.
Then by (4.15.1), 
\begin{equation*}
I' = \{ \dots, a-4, a-4 \}, \quad I = \{ a, a+1, \dots, \}
\end{equation*}
for 
some $a$.  It follows that 
$a$ is the tail of an interval.  One can replace $I$ by 
$J = \{ a, a, \dots\}$ which produces a new symbol $\vL'$
corresponding to $C_{\la',\ve'}$.  This works also for a symbol
$(A, B)$ similar to $\vL$, where $(A, B) \lra E$.   
Now $\la'$ is obtained from $\la$ by replacing two rows of length $h$
by two rows of length $h-1$.  Since $\ve(h-1) = 0$, 
we have $\ve'(h-1) = 0$ and  4.11 can be applied to show that 
$Y_{u,v}^F \ne \emptyset$ for a split element $v \in C_{\la',\ve'}^F$,
and we get the assertion in a similar way as above. 
The other cases are dealt similarly.
\par
Now we may assume that $\la$ consists of even rows.  
Assume that there exists $h$ such that $\ve(h) = 0$, and 
consider a semi-interval $I = I_h^0$.  There are two possibilities 
for $I'$ whose distance is $\le 4$, i.e., $I' = I_{h-2}^1$ or 
$I_{h-2}^0$.  First assume that $I' = I_{h-2}^1$.  Then $I,I'$ is
written as 
\begin{equation*}
I' = \{ \dots, a-6, a-4\}, \quad I = \{ a, a, a+4, a+4, \dots\}
\end{equation*}
by (4.15.1).
Then one can replace the tail $a$ of $I$ by $a-1$, which divide $I$ into
two semi-intervals 
\begin{equation*} 
J_1 = \{ a-1, a\}, \quad J_2 = \{ a+4, a+4, \dots\},
\end{equation*}
and $I'\cup J_1$ form a part of some interval. 
The situation is the same for a symbol 
$(A,B)$ similar to $\vL$.  This produces a new symbol 
$\vL' \lra C_{\la',\ve'}$  where $\la'$ is obtained from $\la$
 by replacing two rows of length $h$ by two rows of length $h-1$
where $\ve'(h-1) = \w$.  Hence  by 4.11, we have 
$Y_{u,v}^F \ne \emptyset$ for a split element $v \in C_{\la',\ve'}$.
Another case is dealt similarly.
\par
Finally we may assume that $\la$ consists of even rows $h$ with 
$\ve(h) = 1$. 
We consider $I = I_h^1$, $I' = I_{h-2}^1$. Then we have
\begin{equation*}
I' =  \{ \dots, a-5, a-3\}, \quad I = \{ a, a+2, a+4, \dots\}.
\end{equation*}
Note that $I, I'$ are a part of a common interval.  By replacing 
$a$ by $a-1$, we have new semi-intervals 
\begin{equation*}
J' = \{ \dots, a-5, a-3, a-1\}, \quad J = \{ a+2, a+4, \dots\}.
\end{equation*}
This produces a new symbol $\vL' \lra C_{\la',\ve'}$, where    
$\la'$ is obtained from $\la$ by replacing one row of length $h$ 
by one row of length $h-2$. Since $\ve(h-2) = 1$, the argument in 
4.11 can be applied, and we get the assertion of the theorem.
Theorem 4.2 is now proved.
\para{4.16.}
We shall prove Theorem 4.3.  So assume that $F$ is of non-split 
type.
First consider the case where $L \ne T$.
In this case the proof is done almost similar to the proof of
Theorem 4.2, by using Lemma 4.13 instead of Proposition 4.10 and 
4.11.
However, we have to be careful for the choice of $v$ (cf. the 
condition of $a$ and $a'$ in 4.12)
in applying Lemma 4.13.  In the case where $p = 2$, this is done
along the line in 4.15, by choosing the decreasing number suitably. 
In the case where $ p \ne 2$,  Proposition 4.10 cannot be applied 
directly, and we have to apply a similar argument as in 4.15.  But
this is easier than the case of $p = 2$; the semi-intervals 
are of the form $I_h^1 = \{ a, a+1, a+2, \dots\}$ or 
$I_h^{\w} = \{ a, a, a+2, a+2, \dots\}$, and only $I_h^1$ gives 
an interval.  The distance of two semi-intervals $I, I'$ is $\ge 2$.  If 
the distance is $\ge 3$, no interaction occurs for $I, I'$ in decreasing
one entry. Then the assertion (i) of the theorem is obtained by 
considering the following $(I,I')$.  
\begin{equation*}
(I_h^{\w}, I_{h-2}^{\w}), \quad (I_h^{\w}, I_{h-1}^1), \quad
(I_h^1, I_{h-1}^{\w}), \quad (I_h^1, I_{h-2}^1). 
\end{equation*} 
The details are omitted.
\par
Next we consider the case where $L = T$.  Hence 
$\CW = W_n'$, and we regard it as a subgroup of 
$W_n = W_n'\lp\s\rp$.
We consider the decomposition of $E = V_{(C',\CE')}$ in 
(1.6.1).  Assume that $E$ is $F$-stable, and let $\wt E$
be the extension of $E$ by $\s_{(C',\CE')}$.  Let 
$\wt E_1$ be the extension of $E_1 = V_{(C_1,\CE_1)}$ through 
$\s_{(C_1, \CE_1)}$. 
Note that in our case $\dim M_{E_1} = 1$, and so $\s_{\CE',\CE}$
is a scalar map.
By multiplying $\s\iv$ on the both
side of (1.6.2), one can write
\begin{equation*}
\s\iv \circ \s_{(C',\CE')}|_{M_{E_1}\otimes E_1}
           = \s_{\CE',\CE_1}\otimes \s\iv\circ\s_{(C_1,\CE_1)}.
\end{equation*}
Since $F$ acts trivially on $A_G(u)$ and $A_M(v)$, the extension 
$\wt{\r\otimes \r_1^*}$ is the trivial extension.  
Thus if $F$ acts trivially on $X_{u',v'}$, then Corollary 1.9 
implies that $\s_{\CE',\CE_1}$ is a scalar map by 
$q^{-d_{C_1,C'} + \dim U_Q}$.
Now assume that $E$ corresponds to the symbol 
$\vL$ and $E_1$ is an irreducible component of $E$ corresponding to 
the symbol $\vL'$, where 
$\vL'$ is obtained from $\vL$ by decreasing an entry by 1, under 
the condition in 4.12. (Note that $\vL$ is not a degenerate symbol
since $E$ is $F$-stable.) Let $\wt E$ and $\wt E_1$ be the preferred 
extensions of $E, E_1$, respectively.  Then it is easy to check 
that $\wt E_1$ occurs in the restriction of $\wt E$ on $W_{n-1}$. 
Thus again by using the arguments in 4.15 (see also the remark 
for the case where $L \ne T$ with non-split case), thanks to 
Lemma 4.13, the verification of Theorem 4.3 (ii) is reduced to 
the case where $n = 2$, namely $G \simeq G_1 \times G_1$, where 
$G_1$ is of type $A_1$, and $F$ acts as a permutation of two 
factors. This case is checked as follows (cf. [S1, Lemma 3.11]).
Since the class $C'$ is $F$-stable, $u'$ is 
the product of two regular elements in $G_1$ or 
the product of two identity elements in $G_1$.
Let $\CB_{u'}$ be the variety of Borel subgroups of $G$ containing
$u'$.   Then $F$ acts trivially on the one dimensional $W_n'$-module
$H^{2d_{u'}}(\CB_{u'})$, where $d_{u'} = \dim \CB_{u'}$.
Hence $W_n$-module $H^{2d_{u'}}(\CB_{u'})$ coincides with the identity 
representation or the long sign representation $\e$ (i.e., 
takes the value $\e(r_{\a}) = 1$, $\e(r_{\b}) = -1$, where 
$r_{\a}$ (resp. $r_{\b}$) is the reflection with respect to the 
short root $\a$ (resp. long root $\b$) of the root system of type
$B_2$.) 
according to the cases where $u'$ is regular or identity. 
The corresponding symbol is $(2; 0)$ for the former, 
and $(12; 01)$ for the latter.   The both are preferred extensions, and 
the theorem follows. 

\par\medskip
\par

\end{document}